\documentclass[a4paper, 10pt, conference]{ieeeconf}      % Use this line for a4 paper

\IEEEoverridecommandlockouts                              % This command is only needed if
\usepackage{amsmath, amsfonts, amssymb }
\usepackage{mathtools}
\usepackage{graphicx}
\usepackage{caption}
\usepackage{booktabs}
\usepackage{multirow}
\usepackage{subcaption}
\usepackage{enumerate}
\usepackage{framed}
\usepackage{tabularx}
\usepackage{color}

\DeclareMathOperator*{\sumsum}{\sum\sum}
\DeclareMathOperator*{\argmin}{arg\min}

\newtheorem{thm}{Theorem}

\newtheorem{lem}[thm]{Lemma}
\def\grad{\nabla}

\def\x{{\boldsymbol{\xi}}}

\def\cC{\mathcal{C}}

\def\cH{\mathcal{H}}

\def\cL{\mathcal{L}}

\def\cN{\mathcal{N}}
\def\cO{\mathcal{O}}

\def\cQ{\mathcal{Q}}

\def\cS{\mathcal{S}}

\def\smskip{\smallskip}

\def\texitem#1{\par\smskip\noindent\hangindent 25pt
               \hbox to 25pt {\hss #1 ~}\ignorespaces}

\def\norm#1{\|#1\|}

\newcommand{\BEAS}{\begin{eqnarray*}}
\newcommand{\EEAS}{\end{eqnarray*}}
\newcommand{\BEA}{\begin{eqnarray}}
\newcommand{\EEA}{\end{eqnarray}}
\newcommand{\BEQ}{\begin{eqnarray}}
\newcommand{\EEQ}{\end{eqnarray}}
\newcommand{\BIT}{\begin{itemize}}
\newcommand{\EIT}{\end{itemize}}
\newcommand{\BNUM}{\begin{enumerate}}
\newcommand{\ENUM}{\end{enumerate}}

\newcommand{\BA}{\begin{array}}
\newcommand{\EA}{\end{array}}

\newcommand{\reals}{\mathbb{R}}
\newcommand{\integers}{\mathbb{Z}}

\newcommand{\dom}{\mathop{\bf dom}}

\newif\ifpagenumbering
\pagenumberingtrue

\pagenumberingfalse

\newsavebox{\theorembox}
\newsavebox{\lemmabox}
\newsavebox{\remarkbox}
\newsavebox{\assbox}
\savebox{\theorembox}{\noindent\bf Theorem}
\savebox{\lemmabox}{\noindent\bf Lemma}
\savebox{\remarkbox}{\noindent\bf Remark}
\savebox{\assbox}{\noindent\bf Assumption}

\def\fprod#1{\left\langle#1\right\rangle}

\title{\LARGE \bf
A Parallel Method for Large Scale Convex Regression Problems
}

\author{NECDET S. AYBAT$^{1}$ and ZI WANG$^{2}$% <-this % stops a space
\thanks{$^{1}$Industrial and Manufacturing Engineering Dept.,
        Penn State University, University Park, PA 16802, USA.
        {\tt\small email: nsa10@psu.edu}}%
\thanks{$^{2}$Industrial and Manufacturing Engineering Dept.,
        Penn State University, University Park, PA 16802, USA.
        {\tt\small email: zxw121@psu.edu}}%
}

\begin{document}

\maketitle
\thispagestyle{empty}
\pagestyle{empty}

%%%%%%%%%%%%%%%%%%%%%%%%%%%%%%%%%%%%%%%%%%%%%%%%%%%%%%%%%%%%%%%%%%%%%%%%%%%%%%%%
\begin{abstract}
Convex regression~(CR) problem deals with fitting a convex function to a finite number of observations. It has many applications in various disciplines, such as statistics, economics, operations research, and electrical engineering. Computing the least squares~(LS) estimator via solving a quadratic program~(QP) is the most common technique to fit a piecewise-linear convex function to the observed data. Since the number of constraints in the QP formulation increases quadratically in $N$, the number of observed data points, computing the LS estimator is not practical using interior point methods when $N$ is very large. The first-order method proposed in this paper carefully manages the memory usage through parallelization, and efficiently solves large-scale instances of CR.% problem with very limited memory.
\end{abstract}

%%%%%%%%%%%%%%%%%%%%%%%%%%%%%%%%%%%%%%%%%%%%%%%%%%%%%%%%%%%%%%%%%%%%%%%%%%%%%%%%
\section{Introduction}
\emph{Convex regression}~(CR) problem is %particularly
concerned with fitting a convex function to a finite number of observations. In particular, suppose that we are given $N$ observations $\{(\pmb{x}_\ell,\bar{y}_\ell)\}_{\ell=1}^N\subset\mathbb{R}^n\times\mathbb{R}$ such that
\begin{align}
\label{data}
\bar{y}_\ell=f_0(\pmb{x}_\ell)+\varepsilon_\ell,\quad \ell=1,\ldots,N,
\end{align}
where $f_0: \mathbb{R}^n \rightarrow \mathbb{R}$ is convex, $\varepsilon_\ell$ is a random variable with $E[\varepsilon_\ell]=0$ for all $\ell$. The objective is to estimate the convex function $f_0$ from the observed data points. %where $\{\epsilon_i\}_{i=1}^N$ is a set of independent and identically distributed~(i.i.d.) random variables.
%our aim is to approximate a convex function $\hat{f}$
CR has many applications in various disciplines, such as statistics, economics, operations research, and electrical engineering. M. Mousavi~\cite{mousavi2013shape} employed %convex regression
CR to estimate the value function under infinite-horizon discounted rewards for Markov chains, which naturally arises in various control problems. In economics field, CR is used for approximating consumers' concave utility functions from empirical data~~\cite{meyer1968consistent}. % which is often concave, as reported in.
Moreover, in queueing network context, for models where the expectation of performance measure is convex in model parameters -see%Chen and Yao
~\cite{chen2001fundamentals}, using Monte Carlo methods to compute the expectation give rise to CR problem~\cite{lim2012consistency}.

The most well-known method for CR is the least squares~(LS) problem,
\begin{align}
\label{infinite_dim}
\hat{f}_N=\argmin_{f\in\cC} \sum\limits_{\ell=1}^{N} \big( f(\pmb{x}_\ell ) - \bar{y}_\ell \big)^2,
\end{align}
where $\cC:=\{f:\mathbb{R}^n\rightarrow\mathbb{R} \hbox{ such that } f \hbox{ is convex}\}$. This infinite dimensional problem %can be equivalently written as
is equivalent to a finite dimensional quadratic problem~(QP),
\begin{align}
\label{original}
& \min_{y_\ell\in\reals,~\pmb{\xi}_\ell\in\reals^n} \sum\limits_{\ell=1}^{N}   \big| y_\ell -\bar{y}_\ell \big|^2\\
\hbox{s.t.}  \quad & {y}_{\ell_1} \geq y_{\ell_2} + {\xi_{\ell_2}}^{T} (\pmb{x}_{\ell_1} - \pmb{x}_{\ell_2}) \quad  1\leq \ell_1 \neq \ell_2\leq N. \nonumber
\end{align}
Indeed, let $\{(y_\ell^*,\xi_\ell^*)\}_{\ell=1}^N$ be an optimal solution to \eqref{original}, it is easy to show that when $N \geq n + 1$, $\{y_\ell^*\}_{\ell=1}^N$ is unique, $\hat{f}_N(\pmb{x}_\ell)=y_\ell^*$ and $\xi_\ell^*\in\partial\hat{f}_N(\pmb{x}_\ell)$ for all $\ell$, where $\partial$ denotes the subdifferential. Moreover, $\hat{f}_N\rightarrow f_0$ almost surely is shown in~\cite{lim2012consistency}; and the convergence rate is established in \cite{groeneboom2001estimation} for one-dimensional case, i.e. $n=1$. LS estimator has some significant advantages over many other estimators proposed in the literature for CR. First, LS estimator is a non-parametric regression method as discussed in \cite{Sen11}, which does not require any tuning parameters and avoids the issue of selecting an appropriate estimation structure. On the other hand, as discussed in~\cite{mousavi2013shape}, methods proposed by Hannah and Dunson~\cite{hannah2011approximate,hannah2013multivariate}, are semi-parametric and require adjusting several parameters before fitting a convex function. Second, LS estimator can be computed by solving the QP in \eqref{original}. Therefore, at least in theory, it can be solved very efficiently using interior point methods (IPM). However, a major drawback of LS estimator in practice is that the number of shape constraints in \eqref{original} is $\mathcal{O}(N^2)$. Consequently, the problem quickly becomes massive even for moderate number of observations: the complexity of each factorization step in IPM is $\mathcal{O}(N^3(n+1)^3)$, and the memory requirement of IPM is $\cO \big( N^2(n+1)^2 \big)$ assuming Cholesky factors are stored - see \cite{Boyd04,Nocedal2006}. %\red{ Additionally, IPM require to store information of , which also causes memory issues in large scale real applications. We have more detailed discussion on memory usage of both IPM and our proposed method in Section~\ref{sec:alcc}. }

In this paper, we develop a methodology for parallel computing the LS estimator on huge-scale CR problems. The proposed method carefully manages the memory usage through parallelization, and efficiently solves large-scale instances of \eqref{original}. %The main contribution of the paper is summarized as following: we proposed a parallel method for solving convex regression problem, particularly for large scale dimension. In the scheme, parallelization serves the purpose of both computation efficiency and memory usage reduction.
Indeed, by regularizing the objective in~\eqref{original}, we ensure the feasibility of primal iterates in the limit, and Lipchitz continuity of gradient of the dual function. These properties lead to the main result, Theorem~\ref{bound}, which provides error bounds on the distance between the LS estimator and the optimal solution to the regularized problem. In the rest of the paper, after examining the dual decomposition for large-scale CR instances, we briefly discuss a first-order augmented Lagrangian method for solving QP subproblems. Finally, we conclude with a number of numerical examples.% in Section~\ref{numerical}.

\section{Methodology}
%In the context of this paper, we
Assume that $\{\varepsilon_\ell\}_{\ell=1}^N$ is uniformly bounded by some $B_\varepsilon>0$, $f_0:\mathbb{R}^n\rightarrow\mathbb{R}\cup\{+\infty\}$ is a convex function, and $\{\pmb{x}_\ell\}_{\ell=1}^N$ is a set of independent and identically distributed (i.i.d.) random vectors in $\mathbb{R}^n$ having a common continuous distribution supported on the $n$-dimensional hypercube $\cH:=[-B_x,B_x]^n\subset \mathbf{ri}\dom(f_0)$ for some $B_x>0$, where $\mathbf{ri}$ denotes the relative interior.
%Since $f_0$ is continuous on $\mathbf{ri}\dom(f_0)$ and $\cH$ is compact, $B_y:=\max_{x\in\cH}f_0(x)<\infty$ exists, and the subgradients of $f_0$ on $\cH$ are uniformly bounded, i.e., $\exists B_\xi>0$ such that $\norm{\xi}_2\leq B_\xi$ for all $\xi\in\partial f_0(x)$ and $x\in\cH$. % We finally assume that the sequence of random variables .

Consider \eqref{original} in the following compact form,
\begin{align}
\label{original_compact}
\min_{\pmb{y}\in\reals^N,~\pmb{\xi}\in\reals^{Nn}} \quad & \tfrac{1}{2} \left\| \pmb{y} - \bar{\pmb{y}} \right\| _2^2\\
\text{s.t. } \qquad & A_1~\pmb{y} +  A_2~\pmb{\xi}\geq 0, \nonumber
\end{align}
where $A_1\in\mathbb{R}^{N(N-1)\times N}$ and $A_2\in\mathbb{R}^{N(N-1)\times Nn}$ are the matrices corresponding to the constraints in \eqref{original}. Let
\begin{align}
\label{tikhonov}
(\pmb{y}^*, \pmb{\xi}^* ) := \argmin\limits_{\pmb{y},~\pmb{\xi}} \left\{   \tfrac{1}{2} \big\| \pmb{y} \big\| _2^2 + \tfrac{1}{2} \big\| \pmb{\xi} \big\| _2^2:\  (\pmb{y}, \pmb{\xi} ) \in \chi^*  \right\},
\end{align}
where $\chi^*$ denotes the set of optimal solutions to \eqref{original_compact}. Let $B_\xi:=\norm{\pmb{\xi}^*}_2$. Moreover, since \eqref{original_compact} is a convex QP, strong duality holds, and an optimal dual solution $\pmb{\theta}^*\in\mathbb{R}^{N(N-1)}$ exists. %to the Lagrangian dual of \eqref{original_compact} exists.
Let $B_\theta>0$ such that $\norm{\pmb{\theta}^*}_\infty\leq B_\theta$ for some optimal dual. The complexity result of the proposed method will be provided in terms of constants $B_\xi$ and $B_\theta$.
\subsection{Separability}
To reduce \textit{curse of dimensionality} and develop a first-order parallel algorithm that can solve~\eqref{original_compact} for large $N$, we use dual decomposition to induce separability. To this aim, we partition $N$ observations into $K$ subsets. Let $\{\cC_i\}_{i=1}^K$ denote the collection of indices such that $|\cC_i|\geq n+1$ for all $i$. To simplify the notation, let $N=K \bar{N}$ for some $\bar{N}>n+1$, $$\cC_i:= \big \{(i-1)\bar{N}+1,~(i-1)\bar{N}+2,\ldots,~i\bar{N} \big\}$$ for $1\leq i\leq K$. Throughout the paper, %for all $1\leq i\leq K$,
$\pmb{y}_i\in\mathbb{R}^{\bar{N}}$ and $\pmb{\xi}_i\in\mathbb{R}^{\bar{N}n}$ denote the sub-vectors of $\pmb{y}\in\mathbb{R}^{N}$ and $\pmb{\xi}\in\mathbb{R}^{Nn}$ corresponding to indices in $\cC_i$, respectively.

For every ordered pair $(\ell_1,\ell_2)$ such that $1\leq\ell_1\neq\ell_2\leq N$, there corresponds a constraint in \eqref{original} represented by a row in the matrices $A_1$ and $A_2$ of formulation \eqref{original_compact}. By dualizing all the constraints in \eqref{original} corresponding to $\ell_1\neq\ell_2$ such that they belong to different sets in the partition, i.e. $\ell_1\in\cC_{i}$, $\ell_2\in\cC_{j}$ and $i\neq j$, we form the partial Lagrangian, \vspace{-2mm}
\begin{align*}
\mathcal{L} \left(\pmb{y}, \pmb{\xi}, \pmb{\theta}\right) :=& \frac{1}{2}\sum\limits_{i=1}^K \norm{\pmb{y}_i - \bar{\pmb{y}}_i}_2^2\\
& - \sumsum\limits_{1\leq i \neq j\leq K} {\pmb{\theta}_{ij}}^{\mathsf{T}} \left( A_1^{ij}
\begin{bmatrix} \pmb{y}_i \\ \pmb{y}_j \end{bmatrix} +A_2^{ij} \begin{bmatrix} \pmb{\xi}_i \\ \pmb{\xi}_j \end{bmatrix}  \right),
\end{align*}
which will lead to the following partial dual function
\begin{align}
\label{subgradient}
g(\pmb{\theta}):=&\min_{\pmb{y}\in\reals^N,~\pmb{\xi}\in\reals^{Nn} } \quad \mathcal{L} \left(\pmb{y}, \pmb{\xi}, \pmb{\theta}\right) \\
\text{s.t.} & \quad  A_1^{ii}~\pmb{y}_i + A_2^{ii}~\pmb{\xi}_i \geq 0,\ i=1,\ldots,K, \nonumber%\norm{\pmb{\xi}_i}_2\leq B_\xi,\; \forall i, \nonumber
\end{align}
where for each $i\in\{1,\ldots,K\}$, $A_1^{ii}$ and $A_2^{ii}$ are formed by rows of $A_1$ and $A_2$, respectively, corresponding to all $(\ell_1,\ell_2)$ such that $\ell_1\neq\ell_2$ and $\ell_1,\ell_2\in\mathcal{C}_i$; similarly, for each $(i,j)$ such that $1\leq i\neq j\leq K$, ${A_1}^{ij} $ and ${A_2}^{ij} $ contain the rows of $A_1$ and $A_2$, respectively, corresponding to $\{(\ell_1,\ell_2):\ \ell_1\in\cC_i,~\ell_2\in\cC_j\}$, %\cup\{(\ell_1,\ell_2):\ \ell_1\in\cC_j,~\ell_2\in\cC_i\}$,
and $\pmb{\theta}_{ij}\in\mathbb{R}^{\bar{N}^2}$ denotes the associated dual variables. % with corresponding constraints in \eqref{original}. Throughout the paper
$\pmb{\theta}$ denotes the vector formed by vertically concatenating $\pmb{\theta}_{ij}$ for $1\leq i\neq j\leq K$.

Note that partial Lagrangian $\cL$ is separable and can be written as $\cL \left(\pmb{y}, \pmb{\xi}, \pmb{\theta}\right)=\sum_{i=1}^K \cL_i\left(\pmb{y}_i, \pmb{\xi}_i, \pmb{\theta}\right)$ for some very simple quadratic functions $\cL_i$. %for $1\leq i\leq K$.
Thanks to the separability of $\cL$, computing the partial dual function $g(\pmb{\theta})$, given in \eqref{subgradient}, is equivalent to solving $K$ quadratic subproblems of the form:
\begin{align}
\label{subgradient_split}
&\min_{\pmb{y}_i\in\reals^{\bar{N}},~\pmb{\xi}_i\in\reals^{\bar{N}n} } \mathcal{L}_i (\pmb{y}_i, \pmb{\xi}_i, \pmb{\theta}) \\
&\qquad \text{s.t.} \quad A_1^{ii}~\pmb{y}_i + A_2^{ii}~\pmb{\xi}_i \geq 0,\nonumber %\norm{\pmb{\xi}_i}_2\leq B_\xi, \nonumber
\end{align}
for $1\leq i\leq K$. %One important property of this relaxation scheme is that given
Given the dual variables $\pmb{\theta}$, since all $K$ subproblems can be computed in parallel, % as one may notice that all subproblems are independent with each other, granted that  are given.
one can take advantage of the computing power of multi-core processors. %with hyper-threading.
In the rest of the paper, we discuss %memory-efficient ways of computing
how to compute a solution to \eqref{original} via solving the dual problem: $\max\{g(\pmb{\theta}):\ \pmb{\theta}\geq 0\}$.

\subsection{Projected Subgradient Method for Dual}
%Various methods can be applied to handle \eqref{subgradient},
One of the most well-known methods for solving the dual problem is the projected subgradient method. %For simplicity let $B_\xi=\infty$ in \eqref{subgradient}, and
Let $\pmb{\theta}^0_{ij}=\pmb{0}$ for all $i$, $j$ such that $i\neq j$. Given the $k$-th dual iterate $\pmb{\theta}^k$, $( \pmb{y}^k, \pmb{\xi}^k )$ denotes an optimal solution to the minimization problem in \eqref{subgradient} when $\pmb{\theta}$ is set to $\pmb{\theta}^k$, %set $\pmb{\theta}=\pmb{\theta}^k$ in \eqref{subgradient}, and let $( \pmb{y}^k, \pmb{\xi}^k )$ be the minimizer of \eqref{subgradient},
and $\pmb{\theta}_{ii}^*$ denotes an optimal dual associated with constraints $A_1^{ii}~\pmb{y}_i + A_2^{ii}~\pmb{\xi}_i \geq 0$ in \eqref{subgradient}. The next dual iterate $\pmb{\theta}^{k+1}$ is computed as follows
\begin{align}
\pmb{\theta}_{ij}^{k+1}  =  \prod_{\cS_{ij}} \left( \pmb{\theta}_{ij}^{k} - t_k\left(  A_1^{ij} \begin{bmatrix} \pmb{y}_i^k \\ \pmb{y}_j^k \end{bmatrix} +A_2^{ij} \begin{bmatrix} \pmb{\xi}_i^k \\ \pmb{\xi}_j^k \end{bmatrix} \right) \right),
\end{align}
where $\Pi_{\cS_{ij}}(.)$ denotes the Euclidean projection on to $$\cS_{ij} = \Big\{ \pmb{\theta}_{ij} \geq \pmb{0}: {\pmb{\theta}_{ij}}^{\mathsf{T}} A_2^{ij} + \begin{bmatrix} {\pmb{\theta}_{ii}^*}^{\mathsf{T}} A_2^{ii}  & {\pmb{\theta}_{jj}^*}^{\mathsf{T}} A_2^{jj} \end{bmatrix}  =\pmb{0} \Big\}.$$
Since $\cL$ is linear in $\pmb{\xi}$, $\dom g$ is non-trivial and is given by the Cartesian product of $\cS_{ij}$'s.
The projected subgradient method is guaranteed to converge in function value with a careful selection of step size sequence $\{ t_k \}_{k=1}^{\infty}$, and it requires $\mathcal{O}( 1/ \epsilon^2)$ iterations to obtain an $\epsilon$-optimal solution~-see~\cite{Nesterov04_1B}. %, or explicitly $k \geq \frac{1}{\epsilon^2} B_3 \sigma_{max}^2\big( \begin{bmatrix} A_3 & A_4 \end{bmatrix} \big) $, where $A_3$ and $A_4$ are defined in \eqref{C}.
%It inflates the computing time so fast that the problem is somewhat impractical to solve, at large scale in particular.
However, due to lack of strong convexity of the objective function in \eqref{original_compact} (not in $\pmb{\xi}$), even if the dual variables converge to an optimal dual solution, the primal feasibility cannot be guaranteed in the limit.
\subsection{Tikhonov Regularization Approach}
In order to ensure feasibility in the limit, which cannot be guaranteed by the subgradient method discussed above, we employ Tikhonov regularization, of which convergence properties were investigated in \cite{engl1989convergence}. Given $\gamma>0$, consider
\begin{align}
\label{regularize}
(\pmb{y}(\gamma),\pmb{\xi}(\gamma))=\argmin_{\pmb{y},~\pmb{\xi}} \quad & \tfrac{1}{2} \left\| \pmb{y} - \bar{\pmb{y}} \right\| _2^2 + \tfrac{\gamma}{2} \left\| \pmb{\xi} \right\| _2^2 \\
\text{s.t. } \quad & A_1\pmb{y} +  A_2{\pmb{\xi}}\geq 0.\nonumber %\norm{\pmb{\xi}_i}_2\leq B_\xi, \forall i.  \nonumber
\end{align}
\vspace{-5mm}

\noindent As $\gamma$ decreases to zero from above, the minimizer $(\pmb{y}(\gamma), \pmb{\xi}(\gamma))$ converges to $(\pmb{y}^*, \pmb{\xi}^* )$ defined in \eqref{tikhonov}.
%\eqref{regularize} with $\gamma =0$.%:= \argmin\limits_{\pmb{y},~\pmb{\xi}}  \Big\{\frac{1}{2} \left\| \pmb{y} - \bar{\pmb{y}} \right\| _2^2:\   A_1\pmb{y} +  A_2{\pmb{\xi}}\geq 0  \Big\}$.
\begin{lem}
\label{continuity}
The minimizer of \eqref{regularize}, $ \pmb{y}(\gamma)$ as a function of regularization parameter $\gamma$, is H\"{o}lder  continuous on $[0,\infty)$,\vspace{-1mm}
\begin{equation}
\label{holder_cond}
\norm{\pmb{y}(\gamma) - \pmb{y}^*}_2 \leq B_\xi \sqrt{\gamma}.
\end{equation}
\end{lem}
\vspace{2mm}
\proof
 Let $\left(\pmb{y}(\gamma),\pmb{\xi}(\gamma)\right)$ be the optimal solution to \eqref{regularize} and $\left(\pmb{y}^*,\pmb{\xi}^*\right)$ be defined as in \eqref{tikhonov}. From the first-order optimality conditions of \eqref{regularize} and \eqref{original_compact}, we have
\begin{equation}
\label{optimality condition1}
\begin{pmatrix} \pmb{y}(\gamma) - \bar{\pmb{y}} \\ \gamma~\pmb{\xi}(\gamma) \end{pmatrix} ^{\mathsf{T}} \begin{pmatrix} \pmb{y}^* - \pmb{y}(\gamma) \\ \pmb{\xi}^* - \pmb{\xi}(\gamma) \end{pmatrix} \geq 0,
\end{equation}
\begin{equation}
\label{optimality condition2}
\begin{pmatrix} \pmb{y}^* - \bar{\pmb{y}} \\ 0 \end{pmatrix} ^{\mathsf{T}} \begin{pmatrix} \pmb{y}(\gamma) - \pmb{y}^*  \\ \pmb{\xi}(\gamma) - \pmb{\xi}^* \end{pmatrix} \geq 0.
\end{equation}
%Recall that $\pmb{\xi}^*=\left[\pmb{\xi}_i^*\right]_{i\in\{1,\ldots,K\}}\in\reals^{Nn}$ and $\pmb{\xi}_i^*=\left[\xi_\ell^*\right]_{\ell\in\cC_i}\in\reals^{\bar{N}n}$.
Note that both $(\pmb{y}(\gamma),\pmb{\xi}(\gamma))$ and $(\pmb{y}^*,\pmb{\xi}^*)$ are feasible to \eqref{original_compact} and \eqref{regularize}. This implies $\norm{\pmb{\xi}(\gamma)}_2\leq \norm{\pmb{\xi}^*}_2$. Summing up \eqref{optimality condition1} and \eqref{optimality condition2}, and using $B_\xi=\norm{\pmb{\xi}^*}_2$, it follows that
%\begin{align*}
%\begin{pmatrix} \pmb{y}(\gamma) -  \pmb{y}^* \\ \gamma~\pmb{\xi}(\gamma) \end{pmatrix} ^{\mathsf{T}} \begin{pmatrix} \pmb{y}^* - \pmb{y}(\gamma) \\ \pmb{\xi}^* - \pmb{\xi}(\gamma) \end{pmatrix} &\geq 0,
%\end{align*}
%or equivalently,
\begin{align*}
\left\|  \pmb{y}(\gamma) -  \pmb{y}^* \right\|_2^2 \leq \gamma \pmb{\xi}(\gamma)^{\mathsf{T}} \big(  \pmb{\xi}^* - \pmb{\xi}(\gamma) \big)\leq\gamma B^2_\xi.
%\gamma \big\| \pmb{\xi}(\gamma) \big\|_2 \big\| \pmb{\xi}^* \big\|_2
%, \\
%( \text{by Cauchy - Schwarz}  )\quad & \leq \quad
%( \text{by boundedness of $\pmb{\xi}$} ) \quad & \leq \quad  N^2 B_\xi^2 \gamma.
\end{align*}
\vspace{-1mm}
\endproof
%This result shows that given an error tolerance $\epsilon>0$, one can accordingly choose the regularization parameter $\gamma\approx \epsilon^2$, more precisely $\gamma= \left(\tfrac{\epsilon}{B_\xi}\right)^2$.
Since the objective function in \eqref{regularize} is strongly convex in both $\pmb{y}$ and $\pmb{\xi}$, Danskin's theorem (see~\cite{Bertsekas99}) implies that the Lagrangian dual function of \eqref{regularize} is differentiable; therefore, one can use gradient type methods to solve the corresponding dual problem. Moreover, strong convexity ensures that, one can solve the primal problem by solving the dual problem. Indeed, let $\pmb{\theta}(\gamma)$ be an optimal solution to the dual problem of \eqref{regularize}, we can recover $(\pmb{y}(\gamma),\pmb{\xi}(\gamma))$ by computing the primal minimizers in \eqref{subgradient} when the dual is set to $\pmb{\theta}(\gamma)$. The discussion above shows that achieving primal feasibility is not an issue provided that we can solve the dual of \eqref{regularize}. This motivates the next section, where we briefly state a first-order algorithm that can efficiently solve the dual of \eqref{regularize}.

\subsection{Accelerated Proximal Gradient~(APG) Algorithm}
Let $\rho:\mathbb{R}^d\rightarrow\mathbb{R}$ be a concave function such that $\grad \rho$ is Lipschitz continuous on $\mathbb{R}^d$ with constant $L$, and $\cQ\subset\mathbb{R}^d$ be a compact convex set. The APG algorithm~\cite{Beck09,Tseng08} displayed in Figure~\ref{fig:apg} is based on Nesterov's accelerated gradient method~\cite{Nesterov04_1B,Nesterov05_1J} and solves $\rho^*=\max\{\rho(\eta):\ \eta\in\cQ\}$. Corollary~3 in~\cite{Tseng08}, and Theorem 4.4 in~\cite{Beck09} show that for all $\ell\geq 1$ the error bound is given by $$0\leq\rho^*-\rho(\eta_\ell)\leq \frac{2L}{(\ell+1)^2} \norm{\eta_0-\eta^*}_2^2,$$ where $\eta_0$ is the initial APG iterate and $\eta^*\in\argmin_{\eta\in\cQ} \rho(\eta)$. Hence, using APG one can compute an $\delta$-optimal solution within at most $\cO(\sqrt{L/\delta})$ APG iterations.
\vspace{-5mm}
\begin{figure}[thpb]
\begin{framed}
{\small
 \textbf{Algorithm APG} \big( $\eta_0$ \big) \\
 Iteration 0: Take $ \eta_0^{(1)}=\eta_1^{(2)}=\eta_0, t_1=1 $\\
 Iteration $\ell$: ($\ell \geq 1$) Compute
 \begin{enumerate}
 \item ${\eta_\ell}^{(1)}=  \Pi_{\cQ}\left( \eta_\ell^{(2)}+\frac{\grad \rho(\eta_\ell^{(2)})}{L} \right)$
 \item $t_{\ell+1}= ( 1+ \sqrt{ 1+ 4t_\ell^2} )/2$
 \item $\eta_{\ell+1}^{(2)}= \eta_\ell^{(1)}+ \frac{t_\ell -1}{t_{\ell+1}} \left(\eta_\ell^{(1)}-\eta_{\ell-1}^{(1)}\right)$
 \end{enumerate}
 }
\end{framed}
\vspace{-2mm}
\caption{Accelerated Proximal Gradient Algorithm}
\label{fig:apg}
\vspace{-3mm}
\end{figure}

In this paper, we will use APG algorithm on a slightly different but equivalent problem to \eqref{regularize}. %With Tikhonov regularization, the primal objective function becomes strongly convex with constant $\sigma=\gamma$, which leads to the differentiability of the dual function. Moreover, the dual function is Lipschitz continuous.Before we state the algorithm, first, let us redefine our problem in favor of the algorithm.
Let $A_3$ and $A_4$ denote the matrices formed by vertically concatenating $A_1^{ij}$ and $A_2^{ij}$, respectively, for $1\leq i\neq j\leq K$; and define
\begin{align}
\label{eq:C}
C=\begin{bmatrix}  A_3 & A_4 \\ I & 0 \end{bmatrix},
\end{align}
where $I\in\mathbb{R}^{N\times N}$ identity matrix. For notational convenience, let $ \pmb{\eta}^{\mathsf{T}}=\begin{bmatrix} \pmb{y}^{\mathsf{T}} & \pmb{\xi}^{\mathsf{T}} \end{bmatrix}$, and consider
\begin{align}
\label{omega_primal}
\min_{\pmb{\eta}\in Q_1} \tfrac{1}{2} \norm{\pmb{y} -  \bar{\pmb{y}}}_2^2 +  \tfrac{\gamma}{2} \left\| \pmb{\xi} \right\| _2^2, \quad \text{s.t.} \quad  C~\pmb{\eta} \geq 0,
\end{align}
where  $Q_1 := \big\{ (\pmb{y},\pmb{\xi}):  A_1^{ii}~\pmb{y}_i +A_2^{ii}~\pmb{\xi}_i \geq  0, 1\leq i \leq K\big\}$. Note that \eqref{omega_primal} is different from \eqref{regularize} only in constraints $\pmb{y}\geq \mathbf{0}$. Via possibly shifting all the observations $\{\pmb{y}_\ell\}_{\ell=1}^N$ up by a sufficiently large quantity, we can assume without loss of generality that $\pmb{y}^*(\gamma)\geq \mathbf{0}$ under bounded error assumption, i.e. $|\varepsilon_\ell|\leq B_\varepsilon$ for all $\ell$. Therefore, \eqref{regularize} and \eqref{omega_primal} are indeed equivalent problems. Consider the dual problem of \eqref{omega_primal},
\begin{equation}
\label{lambda_dual}
 \max_{\pmb{\theta}} g_\gamma(\pmb{\theta}) \quad \hbox{ s.t. } \quad \pmb{\theta} \in Q_2,
\end{equation}
where $Q_2 := \big\{ \pmb{\theta} : \norm{\pmb{\theta}}_2 \leq B_\theta,~\pmb{\theta} \geq 0 \big\}$, and
\begin{align}
\label{eq:g_gamma}
g_\gamma(\pmb{\theta}) = \min\limits_{(\pmb{y}, \pmb{\xi}) \in Q_1 }  \left\{ \tfrac{1}{2} \big\| \pmb{y} -  \bar{\pmb{y}} \big\| _2^2 + \tfrac{\gamma}{2} \left\| \pmb{\xi} \right\| _2^2 - \pmb{\theta}^{\mathsf{T}} C\pmb{\eta} \right\}.
\end{align}
Let $\pmb{\eta}(\pmb{\theta})$ be the minimizer in \eqref{eq:g_gamma}. Theorem 7.1 in~\cite{nesterov2005excessive} and Danskin's theorem imply that
\begin{align}
\label{dual gradient}
\nabla g_\gamma(\pmb{\theta}) = - C \pmb{\eta}(\pmb{\theta})
\end{align}
is Lipschitz continuous with constant
\begin{align}
\label{eq:lischitz}
L_g = \tfrac{1}{\gamma}~\sigma_{\max}^2(C).
\end{align}
Parallel~APG algorithm~(P-APG), displayed in Fig.~\ref{fig:papg}, is the customized version of APG algorithm in Fig.~\ref{fig:apg} to solve \eqref{lambda_dual}. Note that at each iteration computation in Step 1) can be done in parallel using $K$ processors, each solving a smaller QP.
%to compute a $\delta$-optimal solution $\pmb{\theta}_\delta$.%, of Theorem~\ref{bound}.
\vspace{-3mm}
\begin{figure}[thpb]
\begin{framed}
{\small
\textbf{Algorithm P-APG} \big( $\gamma$ \big) \\
Iteration 0: Take $ \pmb{\theta}_0^{(1)} = \pmb{\theta}_1^{(2)} = \pmb{0}, t_1=1 $\\
Iteration $\ell$: ($\ell \geq 1$) Compute
\begin{enumerate}
\item $\pmb{\eta}_\ell = \argmin\limits_{(\pmb{y}, \pmb{\xi}) \in Q_1 }\left\{ \tfrac{1}{2} \big\| \pmb{y} -  \bar{\pmb{y}} \big\| _2^2 + \tfrac{\gamma}{2} \left\| \pmb{\xi} \right\| _2^2 - \left(\pmb{\theta}_\ell^{(2)}\right)^{\mathsf{T}} C\pmb{\eta} \right\}$
\item ${\pmb{\theta}_\ell}^{(1)}=  \Pi_{Q_2}\left( \pmb{\theta}_\ell^{(2)}-\frac{1}{L_g}C\pmb{\eta}_\ell \right)$
\item $t_{\ell+1}= ( 1+ \sqrt{ 1+ 4t_\ell^2} )/2$
\item $\pmb{\theta}_{\ell+1}^{(2)}= \pmb{\theta}_\ell^{(1)}+ \frac{t_\ell -1}{t_{\ell+1}} \left(\pmb{\theta}_\ell^{(1)}-\pmb{\theta}_{\ell-1}^{(1)}\right)$
\end{enumerate}
}
\end{framed}
\vspace{-2mm}
\caption{ Parallel APG Algorithm}
\label{fig:papg}
\vspace{-4mm}
\end{figure}

Note that the iteration complexity of gradient ascent method on \eqref{lambda_dual} is $\cO(L_g/\delta)=\cO(B_\theta^2(\gamma\delta)^{-1})$. On the other hand, P-APG in Fig.~\ref{fig:papg} can compute a $\delta$-optimal solution to \eqref{lambda_dual} within $\cO(\sqrt{L_g/\delta})$ iterations. More precisely, \eqref{eq:lischitz} implies $\cO(B_\theta (\gamma\delta)^{-1/2})$ complexity for P-APG on \eqref{lambda_dual}.
%Note that by using an accelerated method to maximize $g_\gamma$, we enjoy a better iteration complexity than $\cO(\delta^{-2})$ complexity of subgradient method to maximize $g$, while ensuring primal feasibility in the limit.

%Suppose that for given Tikhonov parameter $\gamma$ and APG tolerance parameter $\delta$, we run APG algorithm in Fig. \ref{fig:apg} on \eqref{lambda_dual} to compute a $\delta$-optimal solution $\pmb{\theta}_\delta$. %Let $\pmb{y}^*$ be the unique solution to \eqref{original_compact}, and
Let $\pmb{\theta}_\delta$ be a $\delta$-optimal solution to \eqref{lambda_dual}, and $(\pmb{y}_\delta, \pmb{\xi}_\delta)$ be the optimal solution to the minimization problem in \eqref{subgradient} when $\pmb{\theta}$ is set to $\pmb{\theta}_\delta$. In Theorem~\ref{bound}, which is the main result of this paper, we establish an error bounds on suboptimality $\norm{\pmb{y}_\delta-\pmb{y}^*}_2$, and on infeasibility $\norm{ (A_1 \pmb{y}_{\delta} + A_2 \pmb{\xi}_{\delta} )_{-} }_2$, where $(\pmb{x})_{-}:=\max\{ -\pmb{x},\pmb{0}\}$.   %, \red{ and an error bound of constraint violation,} where.
\begin{thm}
\label{bound}
Let $\left(\pmb{y}(\gamma), \pmb{\xi}(\gamma)\right)$ and $\pmb{\theta}^*$ denote the optimal solutions to \eqref{regularize}  and \eqref{lambda_dual}, respectively. Let $\pmb{\theta}_{\delta}$ be a $\delta$-optimal solution to \eqref{lambda_dual}, and $\left(\pmb{y}_\delta, \pmb{\xi}_\delta\right)$ be the minimizer in %\eqref{subgradient}
\eqref{eq:g_gamma} when $\pmb{\theta}$ is set to $\pmb{\theta}_\delta$. %Then $\pmb{y}_\delta$ is close to the unique optimal solution $\pmb{y}^*$ to \eqref{original_compact}, i.e.
For all $\delta > 0$, the following bounds hold: %the errors in optimality and feasibility are bounded by
\vspace{-2mm}
\begin{flalign}
\label{bound ineq}
&\norm{ \pmb{y}_{\delta} - \pmb{y}^*}_2 \leq B_\xi \sqrt{\gamma} + \sqrt{\tfrac{2\delta}{\gamma}} \sigma_{\max}(C),\\
&\norm{ (A_1 \pmb{y}_{\delta} + A_2 \pmb{\xi}_{\delta} )_{-} }_2 \leq \sqrt{\tfrac{2 \delta }{\gamma}}\sigma_{\max}(C). \label{bound vio}
\end{flalign}
\vspace{-5mm}

%\noindent where $\kappa(C):=\sigma_{\max}(C)/\sigma_{\min}(C)$ is the condition number.% of $C$, i.e. $\kappa(C)$.
\end{thm}
\proof
Since $g_\gamma$ is Lipschitz continuous with constant $L_g$ given in \eqref{eq:lischitz}, we have
\vspace{-3mm}
\begin{align*}
\left\| \nabla g_\gamma( \pmb{\theta}_1) - \nabla g_\gamma( \pmb{\theta}_2) \right\|_2 \leq \frac{\sigma_{\max}^2(C)}{\gamma} \left\| \pmb{\theta}_1 -  \pmb{\theta}_2\right\|_2.
\end{align*}
\vspace{-4mm}

\noindent Moreover, first order optimality conditions for \eqref{lambda_dual} imply
\vspace{-1mm}
\begin{align}
\label{lambda first order}
- \fprod{ \nabla g( \pmb{\theta}^*), \quad \pmb{\theta} - \pmb{\theta}^*} \geq 0, \qquad \forall \pmb{\theta}\in Q_2. %\geq \mathbf{0}.
\end{align}
\vspace{-5mm}

\noindent From (2.1.7) in \cite{Nesterov04_1B}, it follows that
\begin{align*}
 -  \nabla g_\gamma( \pmb{\theta}^*)^{\mathsf{T}} (  \pmb{\theta}_{\delta} - \pmb{\theta}^* ) &+ \frac{\gamma}{2\sigma_{\max}^2(C)}  \left\| \nabla g_\gamma( \pmb{\theta}_{\delta}) - \nabla g_\gamma( \pmb{\theta}^*) \right\|_2^2 \\
& \leq  - g_\gamma( \pmb{\theta}_{\delta}) + g( \pmb{\theta}^*) \leq \delta.
%( \text{by $\delta$-optimal} )& \leq \delta
\end{align*}
Using \eqref{eq:C}, \eqref{dual gradient} and \eqref{lambda first order}, we have
\begin{align}
\left\| \begin{bmatrix} A_3 (\pmb{y}(\gamma) - \pmb{y}_\delta)+A_4(\pmb{\xi}(\gamma) - \pmb{\xi}_\delta) \\ \pmb{y}(\gamma) - \pmb{y}_\delta \end{bmatrix} \right\|_2^2  &\leq \tfrac{2\delta}{\gamma}\sigma_{\max}^2(C). \label{eq:mixed_bound}
\end{align}
Hence, together with \eqref{holder_cond}, it implies \eqref{bound ineq}. Moreover, since $\norm{\pmb{x}-\pmb{y}}_2\geq\norm{(\pmb{x})_{-}-(\pmb{y})_{-}}_2$ for any $\pmb{x}$ and $\pmb{y}$, we also have
%$\norm{C ( \pmb{y}(\gamma) - \pmb{y}_\delta )}_2 \leq \sqrt{\frac{2 \delta }{\gamma}} \sigma_{\max}(C)$.
%from Lemma \ref{LIC}, %and Cauchy - Schwarz Inequality, it follows that
\begin{align*}
 \norm{(A_3\pmb{y}_\delta+A_4\pmb{\xi}_\delta)_{-}-(A_3\pmb{y}(\gamma)+A_4\pmb{\xi}(\gamma))_{-}}_2  \leq \sqrt{\tfrac{2 \delta }{\gamma}} \sigma_{\max}(C).
\end{align*}
Since $(\pmb{y}(\gamma),\pmb{\xi}(\gamma))$ is feasible to \eqref{regularize}, and $(\pmb{y}_\delta,\pmb{\xi}_\delta)\in Q_1$, above inequality implies \eqref{bound vio}.
%Combining \eqref{holder_cond} with \eqref{omega bound}, %in Lemma~\ref{continuity},
%we obtain the result in \eqref{bound ineq}.
%\red{Let $U(\pmb{y}, \pmb{\xi}, \pmb{\theta}) := \frac{1}{2} \big\| \pmb{y} -  \bar{\pmb{y}} \big\| _2^2 + \frac{\gamma}{2} \left\| \pmb{\xi} \right\| _2^2 - \pmb{\theta}^{\mathsf{T}} C\pmb{\eta}$. By strong convexity and optimality of $( \pmb{y}(\gamma), \pmb{\xi}(\gamma), \pmb{\theta}^*)$, we have $\grad U(\pmb{y}(\gamma), \pmb{\xi}(\gamma), \pmb{\theta}^*) = 0$, i.e. $\grad g_{\gamma}(\pmb{\theta}^*) = 0$, and
%\begin{align*}
%\frac{\gamma}{2\sigma_{\max}^2(C)}  \left\| C \begin{bmatrix} \pmb{y}_\delta \\ \pmb{\xi}_\delta \end{bmatrix} \right\|_2^2  &\leq \delta
%\end{align*}
%which implies$ \norm{ A_3 \pmb{y}_{\delta} + A_4 \pmb{\xi}_{\delta} }_2 \leq \sqrt{\frac{2 \delta }{\gamma}} \sigma_{\max}(C)$. Thus, the result in~\eqref{bound vio} follows.}
\endproof
Next, we prove an important technical lemma that will be used later in Theorem~\ref{thm:xi-bound} to show that $\norm{\pmb{\xi}_\delta-\pmb{\xi}^*}_2$ is small. %Moreover, the proof of Lemma~\ref{LIC} helps us understand more on the structure of our problem.
\begin{lem}\label{LIC}
Assuming that $ \{\pmb{x}_i \}_{i=1}^N $ are uniformly sampled at random from set $\phi=\{ \pmb{x} \in \mathbb{R}^n : \| \pmb{x}\|_{\infty} \leq B_x \}$, the matrix $A_4$ in \eqref{eq:C} has linearly independent columns (LIC).
\end{lem}
\proof
Remember $A_4\in\mathbb{R}^{N(N-\bar{N})\times Nn}$ denotes the matrix formed by vertically concatenating all $A_2^{ij}$ for $1\leq i\neq j\leq K$. Note that rows of $[A_1^{ij}~A_2^{ij}]$ correspond to constraints $y_{\ell_1} - y_{\ell_2} + \pmb{\xi}_{\ell_2}^{\mathsf{T}} (\pmb{x}_{\ell_2} - \pmb{x}_{\ell_1}) \geq 0$, where $\ell_1 \in \mathcal{C}_i$, $\ell_2 \in \mathcal{C}_{j}$ and $i\neq j$. For the sake of simplifying the discussion below, without loss of generality, we fix $i=2$ and $j=1$, and focus on the structure of $A_2^{21}$. %By fixing $i,j$, we don't lose any generality.
Let $\hat{A}_2^{21}$ denote the submatrix of $A_2^{21}$ formed by selecting the rows corresponding to $(\ell_1,\ell_2)\in\cC_2\times\cC_1$ such that $\ell_1=\bar{N}+1$ and $\ell_2 \in \mathcal{C}_{1}$. Hence, we have
\begin{align}
\label{A_4}
\hat{A}_2^{21} = \begin{bmatrix} \mathbf{X} & \mathbf{0}\end{bmatrix}
\end{align}
where $\mathbf{0}$ is the matrix of zeros and $\mathbf{X}\in\mathbb{R}^{\bar{N}\times\bar{N}n}$ such that
\begin{align}
\label{X}
\mathbf{X}= \begin{bmatrix}
		 & {\pmb{\bar{x}}_1}^{\mathsf{T}} &\pmb{0}^{\mathsf{T}} &\pmb{0}^{\mathsf{T}} & \ldots & \pmb{0}^{\mathsf{T}} \cr
		 & \pmb{0}^{\mathsf{T}}& {\pmb{\bar{x}}_2}^{\mathsf{T}} &\pmb{0}^{\mathsf{T}} & \ldots & \pmb{0}^{\mathsf{T}}  \cr
		 & \pmb{0}^{\mathsf{T}} & \pmb{0}^{\mathsf{T}} & {\pmb{\bar{x}}_3}^{\mathsf{T}} &\ldots & \pmb{0}^{\mathsf{T}} \cr
		 & \vdots & \vdots & \vdots &\ddots & \vdots \cr
		 & \pmb{0}^{\mathsf{T}} & \ldots & \pmb{0}^{\mathsf{T}} &\ldots & {\pmb{\bar{x}}_{\bar{N}}}^{\mathsf{T}} \cr
\end{bmatrix}
\end{align}
and $\pmb{\bar{x}}_\ell := \pmb{x}_{\ell} - \pmb{x}_{\bar{N}+1}$ for $1\leq\ell\leq\bar{N}$.

Fix $1\leq j\leq K$. Note that for each $\ell\in\cC_j$, there corresponds $n$ columns in $A_4$; and
%Given $1\leq j\leq K$,
the zero structure in \eqref{A_4} implies that each column of $A_4$ corresponding to $\cC_j$ is linearly independent with $\bar{N}n$ columns in $A_4$ corresponding to $\cC_k$ with probability 1~(w.p.~1) for all $k\neq j$. Moreover, when we focus on \eqref{X}, we also see that any one of the $n$ columns in $A_4$ corresponding to $\bar{\ell}\in\cC_j$ is also linearly independent with $n$ columns in $A_4$ corresponding to $\ell\in\cC_j$ with probability 1 for all $\ell\neq \bar{\ell}$. Therefore, to show that $A_4$ has linearly independent columns, it is sufficient to show that for any given $1\leq j\leq K$ and $\ell\in\cC_j$, the corresponding $n$ columns of $A_4$ are linearly independent w.p.~1.

Let $D\in\mathbb{R}^{N(N-\bar{N})\times n}$ be the submatrix of $A_4\in\mathbb{R}^{N(N-\bar{N})\times Nn}$ corresponding to columns $\bar{\ell}\in\cC_j$ for some $1\leq j\leq K$; and $\pmb{d}_{\ell_1 \ell_2}^{\mathsf{T}}$ denote the row of $D$ corresponding to %$1\leq\ell_1, \ell_2\leq N$
$(\ell_1, \ell_2)$ such that $\ell_1$ and $\ell_2$ belong to different sets in the partition. Clearly,
\vspace{-1mm}
\begin{equation}
\pmb{d}_{\ell_1 \ell_2}^{\mathsf{T}}=
\left\{
  \begin{array}{ll}
    \big(\pmb{x}_{\bar{\ell}} - \pmb{x}_{\ell_1} \big)^{\mathsf{T}}, & \hbox{if $\ell_1\not\in\cC_j$ and $\ell_2=\bar{\ell}$;} \\
    \pmb{0}^{\mathsf{T}}, & \hbox{o.w.}
  \end{array}
\right.
\end{equation}
Without loss of generality, we fix $j>1$ and consider $\bar{D} \in \mathbb{R}^{\bar{N} \times n}$ which denotes the submatrix of $D$ corresponding to the rows $\pmb{d}_{\ell_1 \ell_2}^{\mathsf{T}}$ such that $\ell_1 \in \mathcal{C}_1$ and $\ell_2=\bar{\ell}\in \mathcal{C}_j$. The following discussion is true for any $C_k$ such that $k\neq j$, but setting $k=1$ simplifies the notation in $\bar{D}$.
\[
\bar{D} = \begin{pmatrix}  \pmb{x}_{\bar{\ell}}^{\mathsf{T}} - \pmb{x}_1^{\mathsf{T}} \\
						 \vdots \\
						  \pmb{x}_{\bar{\ell}}^{\mathsf{T}} - \pmb{x}_\ell^{\mathsf{T}}\\
						  \vdots \\
						  \pmb{x}_{\bar{\ell}}^{\mathsf{T}} - \pmb{x}_{\bar{N}}^{\mathsf{T}}
\end{pmatrix}
\]
It suffices to show that $\bar{D}$ has LIC. %Consider the linear system $\bar{D}\pmb{\alpha}=\pmb{0}$.
Since $\bar{N}\geq n+1$ and $\{\pmb{x}_\ell\}_{\ell=1}^N$ is a set of i.i.d. random vectors in $\mathbb{R}^n$ having a common \emph{continuous} distribution, it can be shown that there exists $n$ linearly independent rows of $\bar{D}$ w.p.~1. Thus, %we can conclude that
 $A_4$ has LIC.
\endproof
\begin{thm}
\label{thm:xi-bound}
There exists $K_1,K_2>0$ such that \vspace{-1mm}
$$\norm{\pmb{\xi}_\delta-\pmb{\xi}^*}_2\leq K_1\sqrt{\gamma}+K_2\sqrt{\tfrac{\delta}{\gamma}}.$$
\end{thm}
\proof
Since $\pmb{y}^*$ is the unique optimal solution to \eqref{original_compact}, \eqref{tikhonov} implies that $\pmb{\xi}^*=\argmin\{\norm{\pmb{\xi}}_2:\ A_1\pmb{y}^*+A_2\pmb{\xi}\geq \pmb{0}\}$. Similarly, \eqref{regularize} implies that $\pmb{\xi}(\gamma)=\argmin\{\norm{\pmb{\xi}}_2:\ A_1\pmb{y}(\gamma)+A_2\pmb{\xi}\geq \pmb{0}\}$. Hence, for $\pmb{h}(\gamma):=A_1(\pmb{y}^*-\pmb{y}(\gamma))$,
\begin{equation}
\pmb{\xi}(\gamma)=\argmin\{\norm{\pmb{\xi}}_2:\ A_1\pmb{y}^*+A_2\pmb{\xi}\geq \pmb{h}(\gamma)\}.
\end{equation}
Sensitivity of projection onto parametric polyhedral sets was studied in~\cite{Yen95_1J}. Using Theorem~2.1 in~\cite{Yen95_1J} and \eqref{holder_cond}, we have
\begin{equation}
\norm{\pmb{\xi}(\gamma)-\pmb{\xi}^*}_2\leq K\norm{\pmb{h}(\gamma)}_2\leq K\sigma_{\max}(A_1)B_\xi\sqrt{\gamma},
\end{equation}
for some $K>0$. Moreover, \eqref{eq:mixed_bound} and Lemma~\ref{LIC} imply that
\begin{equation*}
\norm{\pmb{\xi}_\delta-\pmb{\xi}(\gamma)}_2\leq\frac{\sqrt{\tfrac{2\delta}{\gamma}}\sigma_{\max}(C)+\sigma_{\max}(A_3)\norm{\pmb{y}(\gamma)-\pmb{y}_\delta}_2}{\sigma_{\min}(A_4)}
\end{equation*}
Hence, $\norm{\pmb{\xi}_\delta-\pmb{\xi}(\gamma)}_2\leq\frac{(\sigma_{\max}(A_3)+1)\sigma_{\max}(C)}{\sigma_{\min}(A_4)}\sqrt{\tfrac{2\delta}{\gamma}}$.
\endproof
%Note that if $A_4$ is a full rank matrix, then $C$ has also full rank. Therefore, Lemma~\ref{LIC} implies that $C$ has linearly independent columns. Since we relax the constraints $C\pmb{\eta}\geq\mathbf{0}$ in \eqref{omega_primal} to form the dual function $g_\gamma$, this key property of $C$ will play a crucial role in proving %the error bound in
%Theorem~\ref{bound}.
\subsection{ALCC - An Augmented Lagrangian Method}
\label{sec:alcc}
Now, we first briefly state a first-order algorithm to directly solve \eqref{regularize}. Let $\bar{B}_y$ and $\bar{B}_\xi$ be given such that $\pmb{y}(\gamma)\in\cQ_y:=\{\pmb{y}:\norm{\pmb{y}-\bar{\pmb{y}}}_2\leq\bar{B}_y\}$, and $\pmb{\xi}(\gamma)\in\cQ_\xi:=\{\pmb{\xi}:\norm{\pmb{\xi}}_2\leq\bar{B}_\xi\}$. Such $\bar{B}_y$ and $\bar{B}_\xi$ can be found easily, if we are given a feasible solution $(\hat{\pmb{y}},\hat{\pmb{\xi}})$, i.e. $A_1\hat{\pmb{y}}+A_2\hat{\pmb{\xi}}\geq\pmb{0}$. Indeed, selecting $\bar{B}_y=\bar{B}$ and $\bar{B}_\xi=B/\sqrt{\gamma}$ works, where $\bar{B}:=\sqrt{\norm{\hat{\pmb{y}}-\bar{\pmb{y}}}_2^2+\gamma\norm{\hat{\pmb{\xi}}}_2^2}$. %, and then we extend to compute projection step in parallel APG algorithm by solving regularized problems of \eqref{subgradient_split}.
ALCC~\cite{Aybat13} computes a solution to \eqref{regularize} by inexactly solving a sequence of subproblems: %of the form
\begin{align}
\label{eq:subproblem}
P_k^*&:=\min\{ P_k(\pmb{y}, \pmb{\xi}):\ \pmb{y}\in\cQ_y, \pmb{\xi}\in\cQ_\xi\},\\
P_k(\pmb{y}, \pmb{\xi})&:= \tfrac{1}{2\mu_k} \left\| \pmb{y}  - \bar{\pmb{y}} \right\| _2^2  + \tfrac{\gamma}{2\mu_k}\| \pmb{\x} \|_2^2 + h_k(\pmb{y}, \pmb{\xi}), \nonumber
\end{align}
where $h_k(\pmb{y}, \pmb{\xi}): = \frac{1}{2} \norm{\left( A_1\pmb{y} + A_2 \pmb{\xi} - \pmb{\theta}_k  \right)_{-}}_2^2$, %$\cQ_y:=\{\pmb{y}:\ \norm{\pmb{y}}_2\leq\sqrt{N}B_y\}$, $\cQ_\xi:=\{\pmb{\xi}: \norm{\pmb{\xi}}_2\leq\sqrt{N}B_\xi\}$,
and $\{\pmb{\theta}_k\}$ sequence is defined in Fig.~\ref{fig:alcc}. %, and $(\pmb{x})_{-}:=\max\{ -\pmb{x},\pmb{0}\}$.
For $k \geq 1$, $h_k( \pmb{y}, \pmb{\xi})$ is convex in $\pmb{y}$ and $\pmb{\xi}$ -see Lemma 3.1 in \cite{Aybat13}. Moreover,
\begin{align*}
\nabla_{\pmb{y}} h_k( \pmb{y}, \pmb{\xi}) &= -{A_1}^{\mathsf{T}} \left( A_1\pmb{y} + A_2 \pmb{\xi} - \pmb{\theta}_k \right)_{-},  \\
\nabla_{\pmb{\xi}} h_k( \pmb{y}, \pmb{\xi}) &= -{A_2}^{\mathsf{T}} \left( A_1\pmb{y} + A_2 \pmb{\xi} - \pmb{\theta}_k  \right)_{-}.
\end{align*}
In addition, $\nabla_{\pmb{y}} h_k( \pmb{y}, \pmb{\xi})$ is Lipschitz continuous in $\pmb{y}$ for all fixed $\pmb{\xi}$ with constant $\sigma_{\max}^2(A_1)$, and $\nabla_{\pmb{\xi}} h_k( \pmb{y}, \pmb{\xi} )$ is Lipschitz continuous in $\pmb{\xi}$  for all fixed $\pmb{y}$ with constant $\sigma_{\max}^2(A_2)$. Hence, $\nabla_{\pmb{y}} P_k( \pmb{y}, \pmb{\xi})$ is Lipschitz continuous in $\pmb{y}$ for all fixed $\pmb{\xi}$ with constant $L_k^y:=\frac{1}{\mu_k}+\sigma_{\max}^2(A_1)$, and $\nabla_{\pmb{\xi}} P_k( \pmb{y}, \pmb{\xi} )$ is Lipschitz continuous in $\pmb{\xi}$  for all fixed $\pmb{y}$ with constant $L_k^\xi:=\frac{\gamma}{\mu_k}+\sigma_{\max}^2(A_2)$.

For given $c>1$ and $\kappa>0$, it is shown in~\cite{Aybat13} that the ALCC algorithm, displayed in Fig.~\ref{fig:alcc}, can compute an $\epsilon$-optimal and $\epsilon$-feasible solution to \eqref{original_compact} within $\cO(\log(\epsilon^{-1}))$ ALCC iterations that require at most $\cO(\epsilon^{-1}\log(\epsilon^{-1}))$ MAPG iterations. The bottleneck step at each MAPG iteration is the matrix-vector multiplication with $A_1 \in \mathbb{R}^{N^2-N \times N}$, $A_2 \in \mathbb{R}^{N^2-N \times Nn}$, ${A_1}^{\mathsf{T}}$ and ${A_2}^{\mathsf{T}}$. Due to specific structures of $A_1$ and $A_2$, without forming $A_1$ and $A_2$ explicitly, we can compute $A_1y$ and ${A_1}^{\mathsf{T}} z$ with $\mathcal{O}(N^2 - N)$ complexity for all $y$ and $z$; $A_2 \xi$ and ${A_2}^{\mathsf{T}} \omega$ with $\mathcal{O} \big( n(N^2 - N) \big)$ for all $\xi$ and $\omega$. Indeed, neither $A_1$ nor $A_2$ is stored in the memory, storing only $\{\pmb{x}_\ell\}_{\ell=1}^N$ is sufficient to be able to compute these matrix-vector multiplications.

\begin{figure}[h!]
\begin{framed}
{\small
 \textbf{Algorithm ALCC} ( $\pmb{y}_0, \pmb{\xi}_0, \mu_1, \tau_1, \alpha^y_1, \alpha^\xi_1$ ) \\
 Iteration 0: Take $\pmb{\theta}^0=\pmb{0}, k=1$\\
 Iteration k: ($k \geq 1$)
 \begin{enumerate}
 \item $L^y_k=\frac{1}{\mu_k}+\sigma_{\max}^2(A_1)$, $L^\xi_k=\frac{\gamma}{\mu_k}+\sigma_{\max}^2(A_2)$
 \item $\ell^{\max}_k= 4\sqrt{\frac{L^y_k \bar{B}_y^2+L^\xi_k\bar{B}_\xi^2}{\tau_k}}$
 \item $ (\pmb{y}_k, \pmb{\xi}_k) = \hbox{\textbf{MAPG}} \big( P_k, L_k^y, L_k^\xi, \pmb{y}_{k-1}, \pmb{\xi}_{k-1}, \alpha^y_k, \alpha^\xi_k, \ell^{\max}_k \big)$
 \item $\pmb{\theta}_{k+1} = \frac{\mu_k}{\mu_{k+1}} \big( A_1\pmb{y}_k + A_2 \pmb{\xi}_k - \pmb{\theta}_k \big)_{-}$
 \item $\mu_{k+1}= c~\mu_k,\ \tau_{k+1}=\tau_k/\left(c~k^{1+\kappa}\right)^2$
 \item $\alpha^y_{k+1}=\alpha^y_{k}/\left(c~k^{1+\kappa}\right)^2,\ \alpha^\xi_{k+1}=\alpha^\xi_{k}/\left(c~k^{1+\kappa}\right)^2$
 \end{enumerate}
}
\end{framed}
\vspace{-2mm}
\caption{Augmented Lagrangian Algorithm ALCC}
\label{fig:alcc}
\vspace{-2mm}
\end{figure}

Note that at each iteration of ALCC in Step 2) MAPG algorithm is called to inexactly solve \eqref{eq:subproblem}. Instead of MAPG, one can also use APG in Fig.~\ref{fig:apg} to inexactly solve \eqref{eq:subproblem}. Within MAPG algorithm, step sizes taken in each block-coordinate are determined by the block Lipschitz constant, i.e. for $y$-coordinate the step size is $1/L^y_k$, while it is $1/L_k^\xi$ for the $\xi$-coordinate. On the other hand, within APG algorithm displayed in Fig.~\ref{fig:apg}, the step sizes taken in each coordinate are equal and determined by the global Lipschitz constant. Thanks to this property of MAPG, we are able to obtain faster convergence in practice in comparison to APG algorithm. When $L_k^y\approx L_k^\xi$, their performance are almost the same; however, when $\max\{L_k^\xi, L_k^y\}/\min\{L_k^\xi, L_k^y\}\gg 1$, since APG uses the global constant L, it takes very tiny steps in one of the block-coordinates.

%For convenience of notation, we define
%\begin{align*}
%h_k( \pmb{y}, \pmb{\xi}, \pmb{\theta} ) :=  \frac{1}{2} {d_{\mathcal{K}}}^2 \left( A_1\pmb{y} + A_2 \pmb{\xi} - \frac{\pmb{\theta}_k}{\mu_k}  \right).
%\end{align*}
\begin{figure}[h!]
\begin{framed}
{\small
 \textbf{Algorithm MAPG} \big( $P,~L^y,~L^\xi,~\pmb{y}_0,~\pmb{\xi}_0,~\alpha^y,~\alpha^\xi,\ell^{\max}$ \big) \\
 Iteration 0: Take $ \pmb{y}_0^{(1)}=\pmb{y}_1^{(2)}=\pmb{y}_0, \pmb{\xi}_0^{(1)}=\pmb{\xi}_1^{(2)}=\pmb{\xi}_0, t_1=1 $\\
 Iteration $\ell$: ($\ell \geq 1$)
 \begin{enumerate}
 \item ${\pmb{y}_\ell}^{(1)}= \Pi_{\cQ_y}\left( \pmb{y}_\ell^{(2)}-\frac{1}{L^y}\nabla_{\pmb{y}} P(\pmb{y}_\ell^{(2)}, \pmb{\xi}_\ell^{(2)})\right)$
\item $ {\pmb{\xi}_\ell}^{(1)}=  \Pi_{\cQ_\xi}\left(\pmb{\xi}_\ell^{(2)}-\frac{1}{L^\xi}\nabla_{\pmb{\xi}} P(\pmb{y}_\ell^{(2)}, \pmb{\xi}_\ell^{(2)})\right)$
 \item $t_{\ell+1}= ( 1+ \sqrt{ 1+ 4t_\ell^2} )/2$
 \item $\pmb{y}_{\ell+1}^{(2)}=  \pmb{y}_\ell^{(1)} + \frac{t_\ell -1}{t_{\ell+1}}  \left(\pmb{y}_\ell^{(1)} - \pmb{y}_{\ell-1}^{(1)}\right)$
 \item $\pmb{\xi}_{\ell+1}^{(2)} = \pmb{\xi}_\ell^{(1)}+ \frac{t_\ell -1}{t_{\ell+1}}\left(\pmb{\xi}_\ell^{(1)}-\pmb{\xi}_{\ell-1}^{(1)}\right)$
 \item \textbf{if} $\norm{\pmb{y}_\ell^{(1)}-\pmb{y}_\ell^{(2)}}_2\leq \alpha^y$ \textbf{and} $\norm{{\pmb{\xi}_\ell}^{(1)}-\pmb{\xi}_\ell^{(2)}}_2\leq \alpha^\xi$
 \item \quad \textbf{return} $\big( \pmb{y}_\ell^{(1)},  \pmb{\xi}_\ell^{(1)}\big)$
 \item \textbf{else if} $\ell=\ell^{\max}$
 \item \quad \textbf{return} $\big( \pmb{y}_\ell^{(1)},  \pmb{\xi}_\ell^{(1)}\big)$
 \item\textbf{end if}% $l=l+1$
 \end{enumerate}
}
\end{framed}
\vspace{-2mm}
\caption{ Modified Accelerated Proximal Gradient Algorithm}
\vspace{-4mm}
\end{figure}
Convergence and rate result of MAPG follow directly from APG in \cite{Beck09} with the help of following lemma.

\begin{lem}
\label{lem:MAPG}
Let $f: \mathbb{R}^m \rightarrow \mathbb{R}$ be a convex function, such that $\grad_{x_1}f( \pmb{x}_1, \pmb{x}_2)$ is Lipschitz continuous with respect to $\pmb{x}_1$ with constant $L_1$, and $\grad_{x_2}f( \pmb{x}_1, \pmb{x}_2)$ is Lipschitz continuous in $\pmb{x}_2$ with constant $L_2$. Then we have
\begin{align*}
f(\pmb{z}_1, &\pmb{z}_2)\leq \: f(\pmb{x}_1, \pmb{x}_2) + L_1 \| \pmb{z}_1 - \pmb{x}_1 \|_2^2 + L_2 \| \pmb{z}_2 - \pmb{x}_2 \|_2^2 \nonumber \\
& + \nabla_{\pmb{x}_1} f(\pmb{x}_1, \pmb{x}_2)^{\mathsf{T}} (\pmb{z}_1 -\pmb{x}_1) + \nabla_{\pmb{x}_2} f(\pmb{x}_1, \pmb{x}_2)^{\mathsf{T}} (\pmb{z}_2 -\pmb{x}_2).
\end{align*}
\end{lem}
\proof
From Lipschitz continuity of $\grad_{x_1}f(.,x_2)$ for each $x_2$ and $\grad_{x_2}f(x_1,.)$ for each $x_1$, it follows that
\begin{align}
\label{lipschitz3}
  f( \pmb{y}_1, \pmb{x}_2) \leq  f( \pmb{x}_1, \pmb{x}_2) &+ \nabla_{\pmb{x}_1} f( \pmb{x}_1, \pmb{x}_2)^{\mathsf{T}} (\pmb{y}_1 -\pmb{x}_1) \nonumber \\
 &+ \tfrac{L_1}{2}  \| \pmb{y}_1 - \pmb{x}_1 \|_2^2, \\
 \label{lipschitz4}
 f( \pmb{x}_1, \pmb{y}_2) \leq  f( \pmb{x}_1, \pmb{x}_2) &+ \nabla_{\pmb{x}_2} f(\pmb{x}_1, \pmb{x}_2)^{\mathsf{T}} (\pmb{y}_2 -\pmb{x}_2) \nonumber \\
 &+ \tfrac{L_2}{2}  \| \pmb{y}_2 - \pmb{x}_2 \|_2^2.
\end{align}
Multiplying \eqref{lipschitz3} and \eqref{lipschitz4} with $\frac{1}{2}$, and summing them up, gives us
\begin{align*}
\tfrac{1}{2} &f( \pmb{y}_1, \pmb{x}_2) + \tfrac{1}{2} f( \pmb{x}_1, \pmb{y}_2) \\
\leq \: &  f( \pmb{x}_1, \pmb{x}_2) + \tfrac{1}{2} \nabla f(\pmb{x}_1, \pmb{x}_2)^{\mathsf{T}} \begin{pmatrix} \pmb{y}_1 -\pmb{x}_1 \\  \pmb{y}_2 -\pmb{x}_2 \end{pmatrix} \\
&+ \frac{L_1}{4} \| \pmb{y}_1 - \pmb{x}_1 \|_2^2 +\frac{L_2}{4} \| \pmb{y}_2 - \pmb{x}_2 \|_2^2.
\end{align*}
Let $\pmb{z}_1 = (\pmb{x}_1+\pmb{y}_1)/2$ and $\pmb{z}_2 = (\pmb{x}_2+\pmb{y}_2)/2$, and by convexity of $f$, we have
\begin{align*}
 f(\pmb{z}_1, \pmb{z}_2) \leq \frac{1}{2} f( \pmb{y}_1, \pmb{x}_2) + \frac{1}{2}&  f( \pmb{x}_1, \pmb{y}_2).
\end{align*}
Combining the last two inequality concludes the proof.
\endproof
Let $\{\pmb{y}^{(1)}_\ell,\pmb{\x}^{(1)}_\ell\}_{\ell\in\integers_+}$ be the iterate sequence generated by MAPG algorithm while running on \eqref{eq:subproblem} starting from $(\pmb{y}_{k-1},\pmb{\xi}_{k-1})$. Using Lemma~\ref{lem:MAPG} and adapting the proof of Theorem 4.4 in~\cite{Beck09}, it can be shown that for all $\ell\geq 1$,
\begin{align*}
0&\leq P_k\left(\pmb{y}^{(1)}_\ell,\pmb{\x}^{(1)}_\ell\right)  - P_k^* \\
 &\leq \frac{4 \left( L_k^y \| \pmb{y}_{k-1} - \pmb{y}_k^*\|_2^2 + L_k^\xi\| \pmb{\x}_{k-1} - \pmb{\x}_k^*\|_2^2 \right)}{(\ell +1)^2},
\end{align*}
where $(\pmb{y}_k^*, \pmb{\x}_k^*)$ is a minimizer of \eqref{eq:subproblem}. Note that we have $\norm{\pmb{y}_{k-1}-\pmb{y}_k^*}_2\leq2\bar{B}_y$ and $\norm{\pmb{\xi}_{k-1}-\pmb{\xi}_k^*}_2\leq2\bar{B}_\xi$. Hence, for all $\ell\geq\ell^{\max}_k$, it is guaranteed that $\left(\pmb{y}^{(1)}_\ell,\pmb{\x}^{(1)}_\ell\right)$ is $\tau_k$-optimal to \eqref{eq:subproblem}.

Note that one can also use ALCC, displayed in Fig.~\ref{fig:alcc}, to compute the primal iterates %$\pmb{\eta}_\ell = (\pmb{y}_\ell, \pmb{\x}_\ell)$ corresponding to $\pmb{\theta}_\ell$
$\pmb{\eta}_\ell$ in Step-1 of P-APG in Fig.~\ref{fig:papg} during the $\ell$-th iteration. In particular, at beginning of every P-APG iteration, $\pmb{\eta}_\ell$ can be computed using ALCC to evaluate $\grad g_\gamma(\pmb{\theta}_\ell^{(2)})$. More importantly, thanks to the separability of regularized \eqref{subgradient_split}, one can do this computation in parallel running ALCC on each one of the $K$ processors.

  Let $N=K\bar{N}$ such that $N\geq n+1$. Below we consider the bottleneck memory requirement for solving \eqref{regularize} in 4 cases: \textbf{a)} P-APG with ALCC computing Step-1 in Fig.~\ref{fig:papg}, \textbf{b)} running ALCC \emph{alone} on \eqref{regularize}, \textbf{c)} P-APG with a primal-dual IPM computing Step-1 in Fig.~\ref{fig:papg}, and \textbf{d)} running IPM \emph{alone} on \eqref{regularize}. The bottleneck for case~\textbf{a)} is determined by Step-2 in Fig.~\ref{fig:papg}, due to dual iterates $\pmb{\theta}$ of size ${(K^2-K)\bar{N}^2+K\bar{N}}$. Similarly, for case \textbf{b)} Step-4 in Fig.~\ref{fig:alcc} requires storing $\pmb{\theta}$ of size $K^2\bar{N}^2$. On the contrary, IPM needs to solve a Newton system in each iteration for both cases \textbf{c)} and \textbf{d)}. Assuming Cholesky factorization is stored, one needs to keep $K$ lower triangular matrices in memory of size $\bar{N}(n+1)$-by-$\bar{N}(n+1)$ for case \textbf{c)}, and to keep 1 lower triangular matrix of size $N(n+1)$-by-$N(n+1)$ for case \textbf{c)}.
  %In short, the information stored in the memory for P-APG method with ALCC is $\cO\big( (K^2- K)\bar{N}^2 \big)$, which is far less than $\cO(K^2 \bar{N}^2 (n + 1)^2) $ for IPM. More detailed memory usage comparison
  Above discussion is summarized in Table~\ref{memory}. Note that running IPM within P-APG reduces the memory requirement significantly by a factor of $K$ in comparison to running IPM alone, e.g. if we partition $N$ observations into $K=10$ subsets and each subproblem requires 1GB of memory, then running IPM alone requires roughly 100GB, while IPM within P-APG requires only 10GB in total.
\begin{table}[htdp]
\centering
\caption{Comparison of Memory Usage}
\begin{tabular}{lll}
\toprule
 & IPM & ALCC \\
 \midrule
Alone & $\cO\big( K^2 \bar{N}^2 (n+1)^2\big)$ & $\cO\big( K^2\bar{N^2} \big)$ \\
P-APG & $\cO\big( K \bar{N}^2 (n+1)^2 + (K^2-K)\bar{N}^2\big)$ & $\cO\big( (K^2-K)\bar{N}^2 \big)$ \\
\bottomrule
\end{tabular}
\label{memory}
\end{table}
\vspace{-4mm}
\section{Numerical Study}
\label{numerical}
In this section, we provide a comparison in Matlab among the following methods: Sedumi, ALCC, Mosek, P-APG with Sedumi, P-APG with ALCC, and P-APG with Mosek, on problem~\eqref{regularize} with increasing dimension. The numerical study is mainly aimed to demonstrate how the performance of each method scales with the dimension of the problem.

\begin{figure}[htbp]
\centering
\includegraphics[width=0.35\textwidth]{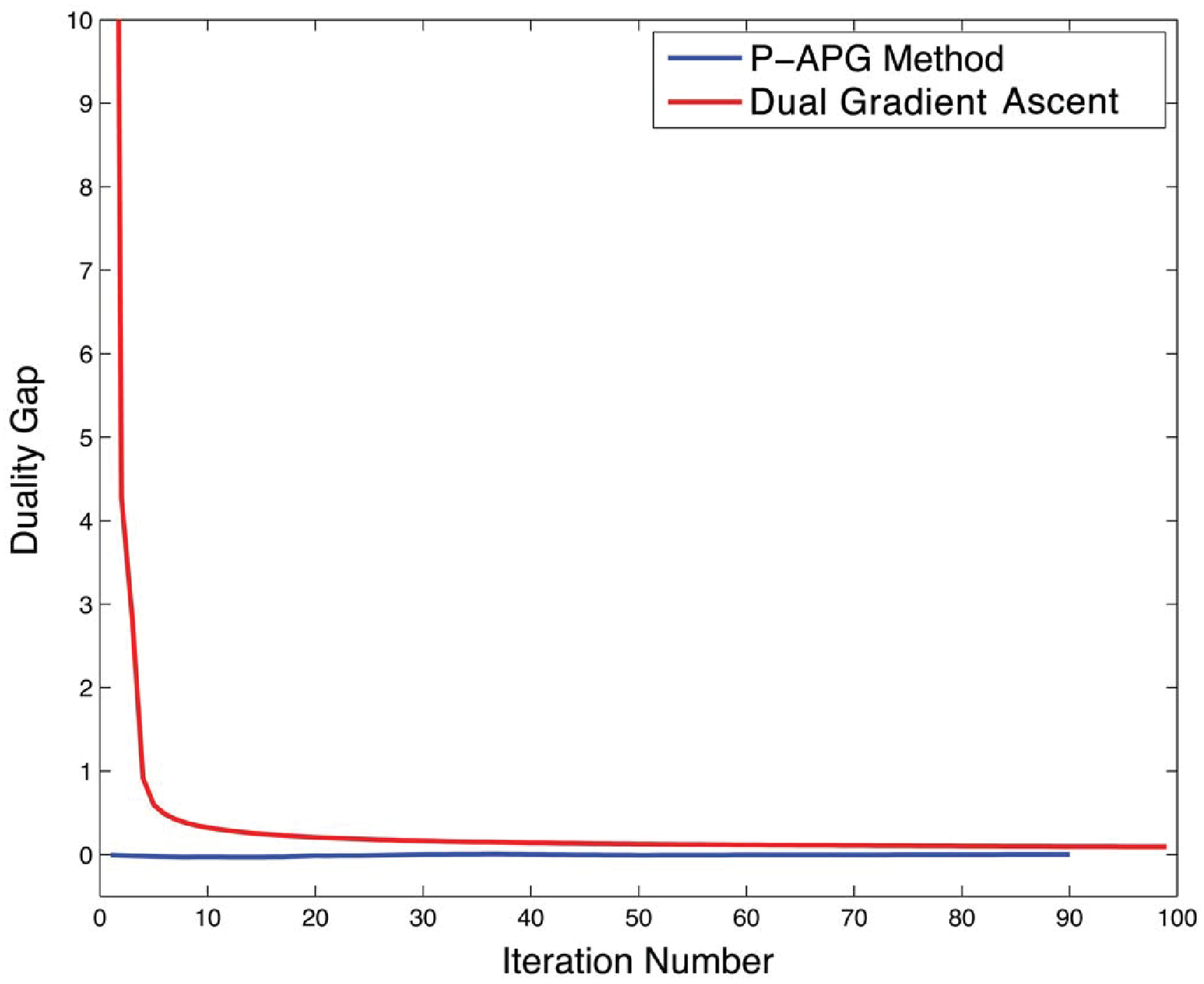}
\caption{ Duality Gap for P-APG and Dual Gradient Ascent}
\label{Gap_All}
\vspace{-4mm}
\end{figure}

\begin{figure}[htbp]
\centering
\includegraphics[width=0.35\textwidth]{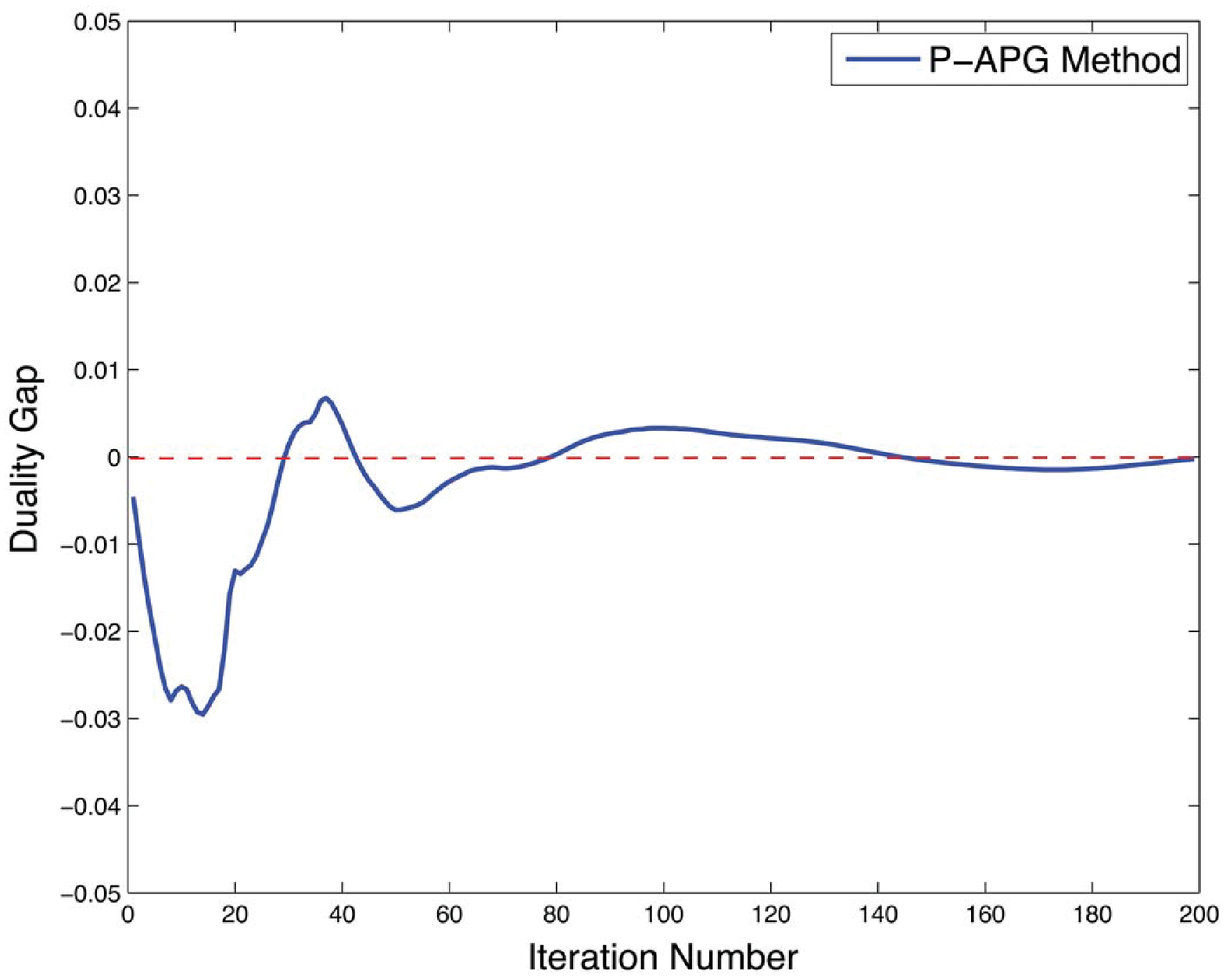}
\vspace{-2mm}
\caption{ Duality Gap for P-APG Method}
\label{Gap_APG}
\vspace{-4mm}
\end{figure}

\begin{figure}[htbp]
\centering
\includegraphics[width=0.35\textwidth]{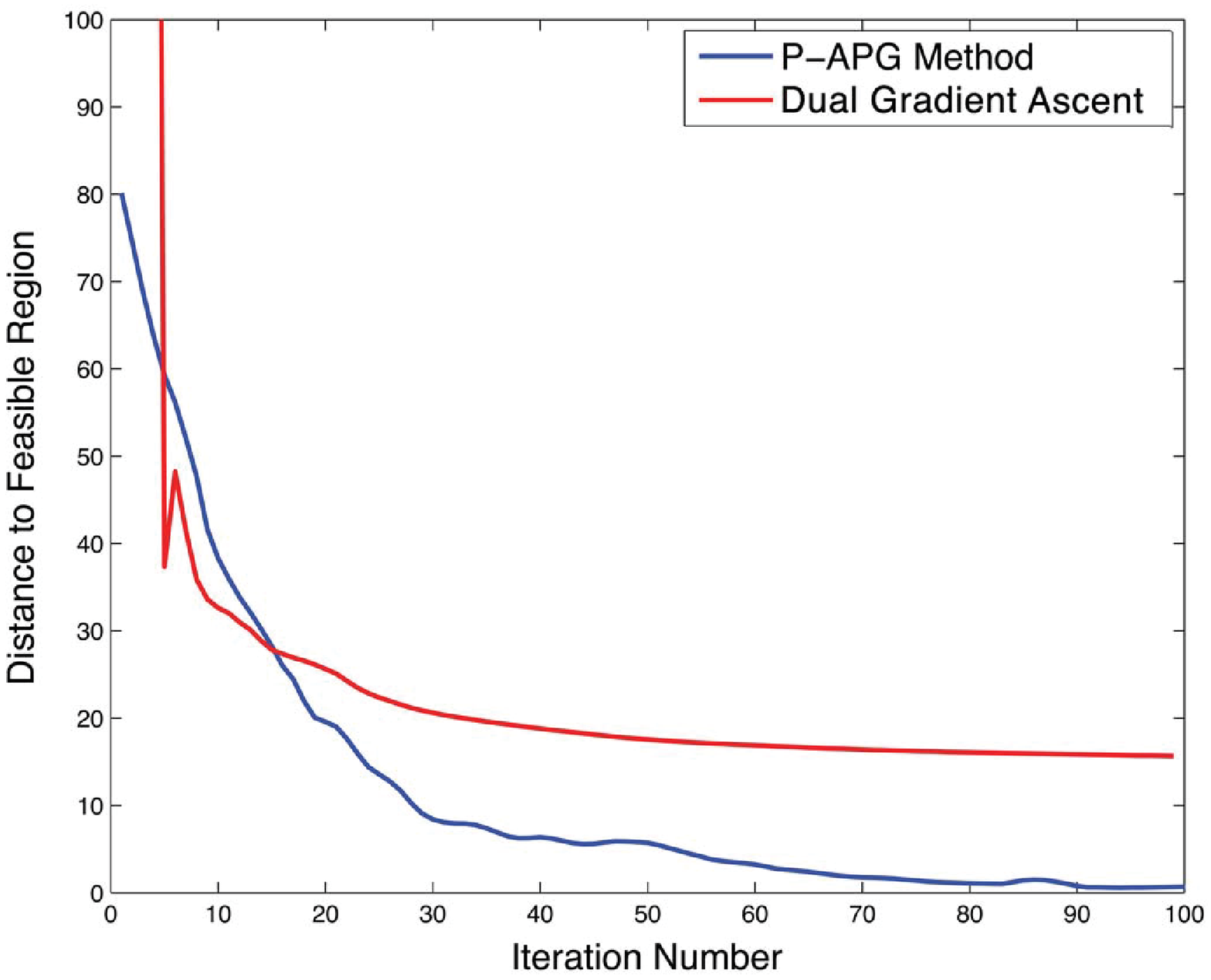}
\caption{ Distance to Feasible Region for P-APG and Dual Gradient Ascent}
\label{Distance_All}
\vspace{-4mm}
\end{figure}

First, we start with a small size problem: $n=5, N=100$. $\Lambda\in\mathbb{R}^{n\times n}$, $\{ \pmb{x}_i \}_{i=1}^N\subset\mathbb{R}^n$ and $\{ \epsilon_i \} _{i=1}^N\subset\mathbb{R}$ are generated randomly with all the components being i.i.d. with $\cN(0,1)$, and $\bar{y}_i$ are generated according to \eqref{data}, where $f_0(\pmb{x})=\frac{1}{2}\pmb{x}^{\mathsf{T}} Q \pmb{x}$, and $Q= \Lambda^{\mathsf{T}} \Lambda$. We %first solve the problem using Sedumi to assess
compare the quality of the solutions computed by P-APG and dual gradient ascent (as the dual function $g_\gamma$ in \eqref{eq:g_gamma} is differentiable). In order to compute dual gradient, $\grad g_\gamma$, one needs to solve $K$ quadratic subproblems. To exploit this parallel structure, we partition the data into two sets, i.e. $K=2$. Within both the dual gradient ascent and P-APG, we called ALCC to compute the dual gradients via solving $K$ QP subproblems. Since we allow violations for the relaxed constraints, ``duality gap" in the paper is defined as $\pmb{\theta}_k^{\mathsf{T}} C \pmb{\eta}_k$ at $k^{th}$ iteration. Fig.~\ref{Gap_All} represents how the duality gap of both methods changes at each iteration. In order to better understand the behavior of P-APG, we report in Fig.~\ref{Gap_APG} the duality gap of P-APG in a smaller scale. Fig.~\ref{Distance_All} reports the infeasibility of iterates, i.e. $\big\| \big( A_1\pmb{y}_k +A_2\pmb{\xi}_k \big)_{-} \big\|_2$. \vspace{-4mm}

\begin{table}[t]
\setlength\tabcolsep{2.7pt}
\caption{ Comparison with test function $\exp( \pmb{p}^{\mathsf{T}}\pmb{x}) $ }
\label{tab:exp}
\centering
{\renewcommand{\arraystretch}{0.65}
\begin{tabular}{llccccc}
\toprule
N & Solver & CPU & W.T. & $\frac{1}{2} \big\| \pmb{y} - \bar{\pmb{y}} \big\|_2^2$ & Gap & Infeas. \\
\midrule
\multirow{6}{*}{200} %& \multirow{2}{*}{Sedumi} & 1.10 & 1.10 & 9.96E-09 & 0 & 0 \\
 & Sedumi & 2.69 & 2.69 & 1.16E-05 & 0 & 0 \\
 & ALCC & 1.86 & 1.86 & 1.16E-05 & -1.19E-07 & 9.6E-02 \\
 & Mosek & 0.74 & 0.74 & 1.16E-05 & 2.68E-08 & 0 \\
 & PAPG(Sedumi) & 29.64 & 14.84 & 1.17E-05 & 1.31E-07 & 9.9E-02 \\
 & PAPG(ALCC) & 9.93 & 4.98 & 1.19E-05 & 2.93E-08 & 9.9E-02 \\
 & PAPG(Mosek) & 3.77 & 1.91 & 1.17E-05 & 9.85E-08 & 9.4E-02 \\
 \midrule
\multirow{6}{*}{400} & Sedumi & O.M. & O.M. & O.M. & O.M. & O.M. \\
 & ALCC & 14.74 & 14.74 & 6.01E-05 & -1.15E-07 & 9.8E-02 \\
 & Mosek & O.M. & O.M. & O.M. & O.M. & O.M. \\
 & PAPG(Sedumi) & 120.48 & 30.28 & 6.00E-05 & 4.52E-09 & 9.7E-02 \\
 & PAPG(ALCC) & 35.83 & 9.11 & 6.11E-05 & -2.87E-08 & 9.8E-02 \\
 & PAPG(Mosek) & 15.73 & 4.12 & 6.00E-05 & 4.97E-09 & 9.7E-02 \\
 \midrule
\multirow{6}{*}{800} & Sedumi & O.M. & O.M. & O.M. & O.M. & O.M. \\
 & ALCC & 93.57 & 93.57 & 2.02E-04 & -8.50E-08 & 9.9E-02 \\
 & Mosek & O.M. & O.M. & O.M. & O.M. & O.M. \\
 & PAPG(Sedumi) & 146 & 19 & 2.46E-04 & 2.41E-08 & 9.9E-02 \\
 & PAPG(ALCC) & 118.54 & 15.77 & 2.10E-04 & 7.01E-08 & 9.8E-02 \\
 & PAPG(Mosek) & 52.43 & 7.52 & 2.05E-04 & 6.49E-08 & 9.7E-02 \\
 \midrule
\multirow{6}{*}{1600} & Sedumi & O.M. & O.M. & O.M. & O.M. & O.M. \\
 & ALCC & N/A & N/A & N/A & N/A & N/A \\
 & Mosek & O.M. & O.M. & O.M. & O.M. & O.M. \\
 & PAPG(Sedumi) & N/A & N/A & N/A & N/A & N/A \\
 & PAPG(ALCC) & 323.68 & 23.94 & 1.97E-03 & 5.85E-09 & 9.9E-02 \\
 & PAPG(Mosek) & 204.56 & 17.14 & 1.85E-03 & -6.72E-10 & 9.9E-02 \\
\bottomrule
\end{tabular} }
\end{table}
\begin{table}[h]
\setlength\tabcolsep{3pt}
\caption{ Comparison with test function $\frac{1}{2}\pmb{x}^{\mathsf{T}}Q\pmb{x}$ }
\label{tab:xqx}
\centering
{\renewcommand{\arraystretch}{0.65}
\begin{tabular}{llccccc}
\toprule
N & Solver & CPU & W.T. & $\frac{1}{2} \big\| \pmb{y} - \bar{\pmb{y}} \big\|_2^2$ & Gap & Infeas. \\
\midrule
\multirow{6}{*}{200} %& \multirow{2}{*}{Sedumi} & 1.10 & 1.10 & 1.68E-09 & 0 & 0 \\
 & Sedumi & 2.88 & 2.88 & 1.25E-04 & 0 & 0 \\
 & ALCC & 3.58 & 3.58 & 1.25E-04 & -3.41E-08 & 9.8E-03 \\
 & Mosek & 0.65 & 0.65 & 1.25E-04 & 2.11E-08 & 0 \\
 & PAPG(Sedumi) & 16.7 & 8.41 & 1.29E-04 & 4.48E-06 & 6.5E-02 \\
 & PAPG(ALCC) & 12.5 & 6.3 & 1.25E-04 & 9.49E-08 & 9.8E-02 \\
 & PAPG(Mosek) & 5.57 & 2.8 & 1.26E-04 & 1.45E-07 & 8.6E-02 \\
 \midrule
\multirow{6}{*}{400} & Sedumi & O.M. & O.M. & O.M. & O.M. & O.M. \\
 & ALCC & 33.3 & 33.3 & 1.02E-03 & -6.84E-08 & 1.4E-02 \\
 & Mosek & O.M. & O.M. & O.M. & O.M. & O.M. \\
 & PAPG(Sedumi) & 164 & 41.2 & 1.02E-03 & -8.62E-08 & 9.5E-02 \\
 & PAPG(ALCC) & 63.2 & 15.7 & 1.01E-03 & 1.53E-07 & 9.7E-02 \\
 & PAPG(Mosek) & 29.01 & 7.43 & 1.02E-03 & -8.74E-08 & 9.5E-02 \\
 \midrule
\multirow{6}{*}{800} & Sedumi & O.M. & O.M. & O.M. & O.M. & O.M. \\
 & ALCC & 140 & 140 & 4.03E-02 & -2.94E-07 & 9.9E-02 \\
 & Mosek & O.M. & O.M. & O.M. & O.M. & O.M. \\
 & PAPG(Sedumi) & 303.33 & 39.87 & 4.04E-03 & -1.95E-07 & 9.9E-02 \\
 & PAPG(ALCC) & 206 & 23.9 & 4.04E-03 & -2.03E-07 & 9.9E-02 \\
 & PAPG(Mosek) & 100.32 & 14.34 & 4.04E-03 & -1.95E-07 & 9.9E-02 \\
 \midrule
\multirow{6}{*}{1600} & Sedumi & O.M. & O.M. & O.M. & O.M. & O.M. \\
 & ALCC & N/A & N/A & N/A & N/A & N/A \\
 & Mosek & O.M. & O.M. & O.M. & O.M. & O.M. \\
 & PAPG(Sedumi) & N/A & N/A & N/A & N/A & N/A \\
 & PAPG(ALCC) & 480 & 29.8 & 5.28E-03 & 1.10E-07 & 9.9E-02 \\
 & PAPG(Mosek) & 273.93 & 21.47 & 5.23E-03 & 5.59E-08 & 9.8E-02 \\
\bottomrule
\end{tabular} }
\end{table}

%By strong convexity of \eqref{regularize}, we can claim that
A primal-dual iterate $(\pmb{\eta},\pmb{\theta})$ is optimal if the duality gap and infeasibility are both zero. As the feasibility happens in the limit, the duality gap in Fig.~\ref{Gap_APG} can go below the red line, which can be explained by the infeasibility of iterates. Therefore, observing a decrease in duality gap only tells one part of the story; without convergence to feasibility, it is not valuable alone as a measure. As shown in the Fig.~\ref{Gap_All} and Fig.~\ref{Gap_APG}, the duality gap converges quickly to zero for both methods. On the other hand, as shown in Fig.~\ref{Distance_All}, constraint violation for P-APG iterates decreases to 0 much faster than it does for the dual gradient ascent iterates. Hence, P-APG iterate sequence converges to the unique optimal solution considerably faster.

The larger scale problems are carried out on a single node at a research computing cluster. The node is composed of one 16-core processor sharing 32GB. For P-AGPG numerical tests, in each job submitted to the computing cluster, an instance of \eqref{regularize} is solved using P-APG on the node such that each subproblem is computed on a different core. The dimension of variables $n=80$ and the number of observations $N=200, 400, 800, 1600$. We partition the set of observations into $K$ subsets. Each one of them consists of 100 points. So, $K=2, 4, 8, 16$ for $N=200, 400, 800, 1600$, respectively. In all the tables, \textit{N/A} means that the wall clock time exceeded 2 hours for the job, and \textit{O.M.} means the algorithm in focus runs out of memory. Also CPU denotes the CPU run time in \emph{minutes}; W.T. stands for wall-clock time in \emph{minutes}. Since the number of constraints increases at the rate of $\mathcal{O}(N^2)$, as the size of problem increases in $N$, we reported the normalized infeasibility and normalized duality gap, which are $\norm{\big( A_1\pmb{y} +A_2\pmb{\xi} \big)_{-}}_2/\sqrt{N^2-N}$ and $\pmb{\theta}_k^{\mathsf{T}} C \pmb{\eta}_k/(N^2-N)$, respectively. We report numerical results for the following test functions: $f_0(\pmb{x})=\frac{1}{2}\pmb{x}^{\mathsf{T}} Q \pmb{x}$, $f_0(\pmb{x}) = \exp( \pmb{p}^{\mathsf{T}} \pmb{x} )$, where $Q$ is generated as discussed before, and $\pmb{p} \in \mathbb{R}^n$ is generated using uniform distribution.

\begin{table}[h!]
\setlength\tabcolsep{3.7pt}
\caption{ Replications with test function  $\exp( \pmb{p}^{\mathsf{T}}\pmb{x}) $ }
\label{tab:exprep}
\centering
{\renewcommand{\arraystretch}{0.55}
\begin{tabular}{lcccccc}
\toprule
 Solver & Rep. & CPU  & W.T.  & $\frac{1}{2} \big\| \pmb{y} - \bar{\pmb{y}} \big\|_2^2$ & Gap & Infeas. \\
\midrule
 \multirow{5}{*}{PAPG(ALCC)} & 1 & 118.54 & 15.77 & 2.10E-04 & 7.01E-08 & 9.8E-02 \\
  & 2 & 131.93 & 17.56 & 3.90E-04 & 5.52E-09 & 9.9E-02 \\
  & 3 & 136.50 & 18.01 & 4.69E-04 & -3.87E-09 & 9.8E-02 \\
  & 4 & 126.43 & 16.74 & 2.91E-04 & -2.53E-09 & 9.9E-02 \\
  & 5 & 144.31 & 18.98 & 5.35E-04 & -7.62E-08 & 9.8E-02 \\
 \midrule
 \multirow{5}{*}{PAPG(Mosek)} & 1 & 52.43 & 7.52 & 2.05E-04 & 6.49E-08 & 9.7E-02 \\
  & 2 & 57.33 & 8.12 & 3.77E-04 & 9.39E-09 & 9.9E-02 \\
  & 3 & 61.53 & 8.64 & 4.54E-04 & -5.30E-09 & 9.9E-02 \\
 & 4 & 55.86 & 8.00 & 2.84E-04 & -7.48E-09 & 9.9E-02 \\
 & 5 & 65.04 & 9.15 & 5.15E-04 & -7.74E-08 & 9.8E-02 \\
\bottomrule
\end{tabular} }
\end{table}
\begin{table}[h!]
\setlength\tabcolsep{3.7pt}
\caption{ Replications with test function $\frac{1}{2}\pmb{x}^{\mathsf{T}}Q\pmb{x}$ }
\label{tab:xqxrep}
\centering
{\renewcommand{\arraystretch}{0.55}
\begin{tabular}{lcccccc}
\toprule
 Solver & Rep. & CPU  & W.T.  & $\frac{1}{2} \big\| \pmb{y} - \bar{\pmb{y}} \big\|_2^2$ & Gap & Infeas. \\
\midrule
 \multirow{5}{*}{PAPG(ALCC)} & 1 & 206.00 & 23.90 & 4.04E-03 & -2.03E-07 & 9.9E-02 \\
  & 2 & 213.27 & 27.63 & 1.00E-03 & -9.04E-08 & 9.6E-02 \\
  & 3 & 211.18 & 27.37 & 1.11E-03 & -1.09E-07 & 9.9E-02 \\
  & 4 & 178.77 & 23.41 & 7.29E-04 & -6.20E-08 & 9.9E-02 \\
  & 5 & 200.62 & 26.14 & 1.27E-03 & -1.05E-07 & 9.6E-02 \\
 \midrule
 \multirow{5}{*}{PAPG(Mosek)} & 1 & 100.32 & 14.34 & 4.04E-03 & -1.95E-05 & 9.9E-02 \\
 & 2 & 79.27 & 10.87 & 9.88E-04 & -1.04E-07 & 9.9E-02 \\
 & 3 & 83.11 & 11.43 & 1.09E-03 & -1.19E-07 & 9.9E-02 \\
 & 4 & 68.90 & 9.66 & 7.16E-04 & -2.30E-08 & 9.9E-02 \\
 & 5 & 79.55 & 10.99 & 1.24E-03 & -1.19E-07 & 9.8E-02 \\
\bottomrule
\end{tabular} }
\end{table}
\addtolength{\textheight}{- 12.3 cm}

%For reference, we also run Sedumi on problem \eqref{original}. As shown in Table~\ref{tab:exp} and~\ref{tab:xqx}, the first and second row of Sedumi with dimension $n=80$ and $N=200$ are the results for~\eqref{original} and~\eqref{regularize}, respectively.
All the algorithms are terminated either when they compute an iterate with normalized infeasibility and normalized duality gap are less than 1E-01 and 1E-06, respectively, or at the end of 2 hours. The numerical results reported in Table~\ref{tab:exp} and~\ref{tab:xqx} show that P-APG solution is very close to the real optimal solution of \eqref{regularize}. Note that ALCC fails to terminate within in 2 hours when $N=1600$; and interior point methods fail to run anything beyond $N=200$ due $\cO(N^2n^2)$ memory requirement. Moreover, in order to test the robustness of P-APG, we solved 5 random instances when $N=800$, of which results are reported in Table~\ref{tab:exprep} and Table~\ref{tab:xqxrep}. %We observe that the run time for both test functions are very close among all 5 replications. As the results of different $N$ values have the same trend, we only provide the results for $N=800$ due to the limited space.
Numerical results show that advantages of P-APG over running IPM or ALCC alone on \eqref{regularize} become more and more evident as the dimension of the problem increases. %: it requires no tuning on the parameters; it is more stable on a variety of problems -- converging to the optimal solution regardless of the dimension of the problem; and most importantly, it is memory efficient, which
\section{Conclusion}
In this paper, we proposed P-APG method to efficiently compute the least squares estimator for large scale convex regression problems. By relaxing constraints partially, we obtained the separability on the corresponding Lagrangian dual problem. Using Tikhonov regularization, we ensured the feasibility of iterates in the limit, and we provided error bounds on 1) the distance between the inexact solution to the regularized problem and the optimal solution to the original problem, 2) the constraint violation of the regularized solution. The comparison in the numerical section demonstrates the efficiency of P-APG method on memory usage compared to IPM. Furthermore, the extended random tests show the stability of P-APG method. Due to limited space, we could not include computational results on real-life data; but they will be made available online at authors' webpage.

%\textbf{Estimation of Option Pricing Function}
%We demonstrate an application of our proposed method to estimate a convex function for European call option pricing. The call option pricing function for European style option is a function of underlying asset price $S_t$ at time $t$, the strike price $x$, the time to expiration $\tau$, the deterministic risk free interest rate $r_{t, \tau}$ and the correspondent dividend yield $\delta_{t, \tau} $ of the asset, which is given by
%\begin{align*}
%C(x)= e^{-\tau r_{t, \tau}} \int_{0}^{\infty} \max(S_T - x, 0) p^*(S_T | S_t, \tau, r_{t, \tau}, \delta_{t, \tau} ) d S_T
%\end{align*}
%Yacine et. al. \cite{ait2003nonparametric} showed that $C$ must be a decreasing function and convex on $x$ in order to rule out arbitrages. We study the real data of NASDAQ-100 call option provided by Yahoo Finance. Figure\ref{option} shows the estimated pricing function with respect to $x$.
%
%\begin{figure}[htbp]
%\centering
%\includegraphics[width=0.5\textwidth]{option}
%\caption{Convex Estimation of NASDAQ-100 Call Option}
%\label{option}
%\end{figure}

%%%%%%%%%%%%%%%%%%%%%%%%%%%%%%%%%%%%%%%%%%%%%%%%%%%%%%%%%%%%%%%%%%%%%%%%%%%%%%%%

%References are important to the reader; therefore, each citation must be complete and correct. If at all possible, references should be commonly available publications.
\bibliographystyle{unsrt}
\bibliography{paper}

\end{document}